%% file: ex_article.tex
\documentclass[hidelinks,onefignum,onetabnum]{siamart250211}

\input{ex_shared}

\ifpdf
\hypersetup{
  pdftitle={LOR Preconditioning for \spectralhp},
  pdfauthor={P. Khurana, H. W\"ustenberg, D. Moxey, A. Chatzopoulos, J. Hoessler, and S. Sherwin}
}
\fi

\begin{document}

\maketitle

\begin{abstract}
Low-order refined (LOR) preconditioning replaces a high-order operator with a low-order discretisation on a refined nodal mesh. For tensor-product elements, the two operators are spectrally equivalent with bounds independent of the polynomial order $P$, but the construction does not extend directly to simplex and mixed-element discretisations. 
This work makes two contributions: it extends LOR preconditioning to simplex and mixed-element discretisations, including triangular, tetrahedral, and prismatic elements, and establishes a generalised Vandermonde transformation linking the modal and nodal LOR formulations, showing that the resulting preconditioned spectra and Krylov convergence are independent of the high-order basis.
Numerical experiments show controlled iteration growth on triangular meshes despite increasing condition number, and controlled iteration counts up to \(P=5\) on tetrahedral, prismatic, and mixed-element meshes. A single algebraic multigrid V-cycle per outer iteration gives the best balance of iteration count and cost. The method is applied to a production incompressible Navier-Stokes simulation of a race-car front-wing and wheel configuration, discretised on a mesh of \(2.87\times10^6\) mixed prismatic and tetrahedral elements giving \(32.2\times10^6\) pressure degrees of freedom at \(P=3\). LOR reduces the mean pressure conjugate gradient (CG) iteration count from \(235.3\) to \(5.5\), and the pressure-solve time over 1000 timesteps by 16.1\%, relative to the default production static-condensation diagonal preconditioner in \textit{Nektar++}. The one-time cost of constructing the LOR preconditioner is amortised over the production simulation.
\end{abstract}

\begin{keywords}
\spectralhp methods, Preconditioning, Low-order refined, Pressure Poisson, Incompressible Navier-Stokes equations
\end{keywords}

\begin{MSCcodes}
65F08, 65N30, 65N35
\end{MSCcodes}

\section{Introduction}
\label{sec:introduction}
\input{main_introduction}

The paper is organised as follows. The LOR methodology is described in \cref{sec:methodology}. The algorithmic formulation and associated design choices are presented in \cref{sec:alg}. The iterative behaviour and conditioning properties of the proposed approach are examined in \cref{sec:LOR_conditioning_iterative}. Scaling performance, hardware efficiency, and runtime characteristics are assessed in \cref{sec:scaling_computational}. Finally, conclusions and directions for future work are given in \cref{sec:conclusions}.

\section{Methodology}
\label{sec:methodology}
\input{main_methodology}

\section{Algorithm}
\label{sec:alg}
\input{algorithm}

\section{Conditioning and Iterative Performance}
\label{sec:LOR_conditioning_iterative}
\input{main_conditioning}
\input{main_iterative}

\section{Scalability and Runtime Performance}
\label{sec:scaling_computational}
\input{main_scaling}

\input{main_computational}

\section{Conclusions}
\label{sec:conclusions}

This work extends low-order refined (LOR) preconditioning to the simplicial and mixed-element discretisations used by \nekpp and establishes its independence from the choice of high-order basis. Existing LOR formulations are typically posed for nodal tensor-product bases, for which the high- and low-order operators are spectrally equivalent with bounds independent of $P$. For a modal hierarchical basis, we show that the LOR preconditioner is related to its nodal counterpart by a congruence transformation under the generalised Vandermonde matrix $\mathcal{V}$. When the high-order operator is transformed consistently, the modal and nodal preconditioned operators are related by a similarity transformation and therefore have identical spectra. The corresponding preconditioned CG iterations are algebraically equivalent under the basis transformation in exact arithmetic. This result is independent of element type and therefore also applies where no spectral-equivalence theory is available. We further develop simplicial subdivision procedures to construct LOR meshes for triangular, tetrahedral, and prismatic elements from the interpolation points of the high-order element. For these element types, no corresponding spectral-equivalence bound is currently available, so their conditioning and iterative behaviour are assessed numerically over the polynomial orders considered here.

We assess the conditioning and iterative performance of the LOR preconditioner on canonical two- and three-dimensional geometries, including mixed-element meshes. Earlier work on quadrilateral and curved-quadrilateral meshes established asymptotically optimal iteration behaviour and an increasing advantage over diagonal preconditioning under mesh refinement \cite{khurana2023comparison}. On triangular meshes, the condition number increases with \(P\), but iteration growth remains controlled. In three dimensions, iteration counts remain controlled up to \(P=6\) on tetrahedral, prismatic, and mixed-element meshes, with LOR consistently requiring fewer iterations than the default \nekpp configuration. For the low-order solve, a single algebraic multigrid V-cycle per outer iteration is used as a practical inexact solve in the tests reported here.

The method is further evaluated on industrial race-car aerodynamic configurations through iteration counts, wall-clock time, and parallel scalability. A single AMG V-cycle per outer iteration was used throughout these tests, having provided a favourable balance between iteration count and computational cost. On the IFW-W geometry, consisting of \(2.87\times10^{6}\) mixed prismatic and tetrahedral elements, LOR reduces the mean pressure CG iteration count at \(P=3\) from \(235.3\) to \(5.5\), with the observed range decreasing from \(32\)--\(5156\) to \(4\)--\(74\). Over \(1000\) timesteps of a production incompressible Navier--Stokes simulation, the pressure-solve time decreases by \(16.1\%\), relative to the default preconditioner. Strong-scaling experiments on ARCHER2 show that this iterative advantage persists with increasing parallelism: outer iteration counts remain essentially independent of MPI rank count, varying from \(43\)--\(45\) at \(P=2\) and \(54\)--\(55\) at \(P=3\). The pressure solve continues to benefit from additional parallelism over most of the rank range considered, with diminishing returns at the highest core counts as communication costs increasingly limit scaling. Although construction of the AMG hierarchy becomes more expensive as the rank count increases, this cost is incurred once and is amortised over the thousands of timesteps required in a production simulation. Together, the strong-scaling and production results show that the reduction in iteration count translates into lower pressure-solve time at industrial scale.

The increase in condition number on simplices is associated with the geometry of the low-order mesh. As \(P\) increases, the interpolation points cluster towards simplex vertices and edges, producing increasingly anisotropic sub-elements. Alternative or optimisation-based node distributions \cite{Chalmers2018}, provide a direct way to address this source of degradation and warrant further study at higher polynomial orders. The present implementation also requires an assembled low-order operator because the multigrid hierarchy is constructed from its matrix entries. Partial-assembly approaches could remove this requirement and enable matrix-free LOR preconditioning \cite{PAZNER2020, Pazner2022GPU}. Porting the implementation to the ongoing redesign of \nekpp would enable GPU execution and the use of GPU-capable multigrid libraries such as MueLu \cite{MueLuURL} and AMGx \cite{AmgX2015}; a direct Hypre interface would also permit finer control of the AMG configuration. The LOR preconditioner is already available for the velocity Helmholtz systems in the solver, but their performance is outside the scope of the present study, which focuses on the pressure Poisson system.

\section*{Acknowledgments}
This project received funding from the European Union’s Horizon 2020 research and innovation programme under the Marie Skłodowska-Curie grant agreement No 955923. The authors gratefully acknowledge Hari Sundar for his expert guidance on linear solvers and multigrid methods, which contributed significantly to the methodology. The authors acknowledge computational resources and support provided by the Imperial College Research Computing Service (\url{http://doi.org/10.14469/hpc/2232}). This work also used the ARCHER2 UK National Supercomputing Service (\url{https://www.archer2.ac.uk}) via the UK Turbulence Consortium (Grant no. EP/R029326/1).

\appendix
\input{appendix}

\bibliographystyle{siamplain}
\bibliography{references}
\end{document}

%% file: ex_shared.tex
\usepackage{amsfonts}
\usepackage{graphicx}
\usepackage{epstopdf}
\usepackage{algorithmic}
\ifpdf
	\DeclareGraphicsExtensions{.eps,.pdf,.png,.jpg}
\else
	\DeclareGraphicsExtensions{.eps}
\fi

\input{includes}
\usepackage{tikz}
\usepackage{enumitem}
\usepackage{needspace}
\usepackage{float}
\usepackage{subcaption}
\usepackage{dsfont}
\usepackage{multirow}
\usepackage{xcolor}
\usepackage{listings} 
\usepackage{lstautogobble}
\usepackage{xspace}
\usetikzlibrary{calc}

\newsiamremark{remark}{Remark}
\newsiamremark{hypothesis}{Hypothesis}
\crefname{hypothesis}{Hypothesis}{Hypotheses}
\newsiamthm{claim}{Claim}
\newsiamremark{fact}{Fact}
\crefname{fact}{Fact}{Facts}

\headers{LOR Preconditioning for SPECTRAL/HP ELEMENT METHOD}{P. Khurana et al.}

\title{LOW-ORDER REFINED PRECONDITIONING FOR SPECTRAL/HP ELEMENT METHOD FOR COMPLEX, 3D GEOMETRIES\thanks{
		\funding{This project received funding from the European Union’s Horizon 2020 research and innovation programme under the Marie Skłodowska-Curie grant agreement No 955923.}}}

\author{Parv Khurana \footnotemark[5] \footnotemark[6] \thanks{corresponding author, \email{Parv621@gmail.com}.}
\and Henrik W\"ustenberg \footnotemark[6] \thanks{corresponding author, \email{h.wustenberg@imperial.ac.uk}.}
\and David Moxey\thanks{Department of Engineering, King’s College London, WC2R 2LS London.}
\and Athanasios Chatzopoulos\footnotemark[5]
\and Julien Hoessler\thanks{CFD Methodology Group, McLaren Racing Limited, GU21 4YH Woking.}
\and Spencer Sherwin\thanks{Department of Aeronautics, Imperial College London, SW7 2AZ London.}}

\usepackage{amsopn}


%% file: includes.tex
\newcommand{\nekpp}{{\em Nektar++}\xspace}

\newcommand{\spectralhp}{spectral/$hp$ element\xspace}
\newcommand{\GLL}{Gauss-Lobatto-Legendre\xspace}
\newcommand{\IncNS}{incompressible Navier-Stokes\xspace}

%% file: main_introduction.tex
\label{sec:methods}
The iterative solution of the \spectralhp method \cite{karniadakis2005spectral} is more computationally expensive than its lower-order counterparts. Firstly, using the \spectralhp discretisation results in denser matrix systems. Higher-order (HO) problems exhibit increased coupling in the degrees of freedom (DOF) as the polynomial order ($P$) increases. This increased coupling implies a higher number of non-zero entries per row in the matrix system, leading to a denser, more complex system. Moreover, the diverse length scales in the geometry and the anisotropic nature of the unstructured meshes needed for complex, three-dimensional, industrial geometries contribute to ill-conditioned matrices. Solving these large, dense, anisotropic, ill-conditioned linear systems becomes difficult and time-consuming for iterative solvers. A new preconditioning technique is therefore required that both accelerates iterative convergence and scales efficiently on modern multi-core HPC systems. 

We address this by developing a \textit{Low-order Refined} (LOR) preconditioner within the \IncNS equations solver of the open-source \spectralhp method framework \nekpp \cite{Moxey2020Nektar}. The LOR preconditioner involves using a corresponding low-order (\textit{P}=1) ``\textit{refined}", spectrally equivalent \cite{CanutoSpectral2007}  finite element discretisation, which is often a structured grid of \GLL (GLL) nodal points. The motivation for this approach is to obtain a sparser preconditioner by leveraging the sparsity of low-order operators and stiffness matrices relative to their higher-order counterparts. The advantages of using a low-order discretisation to improve conditioning properties for the finite element method were observed previously in the works \cite{deville1985chebyshev, deville1990finite, deville1992fourier}. A bounded iterative condition number behaviour with increasing problem size was later established, see references \cite{CANUTO1985315, KimPatera1996, KimPatera1997}, with the bound being $\approx \pi^2/4$ under specific conditions. Another advantage of using a LOR discretisation is the constant memory requirement per degree of freedom, which makes it less memory-intensive than higher-order discretisations \cite{PAZNER2020}. The low memory requirement can be particularly beneficial when dealing with large-scale problems or when computational resources are limited. These methods exhibit minimal sensitivity to high-aspect-ratio elements \cite{PedroFischer2019}. 

Recent years have witnessed significant progress in the development of such low-order preconditioning techniques, especially for meshes with rectangular and hexahedral elements with Lagrange polynomial basis and GLL-based quadrature \cite{Olson2007, Canuto2009, PedroFischer2019, Pazner2022deRham, Phillips2021GPU, Pazner2022GPU}. The performance and conditioning of $\mathds{P}1$ and $\mathds{Q}1$ discretisations in two and three dimensions have been analysed \cite{Canuto2009}, and $\mathds{P}1$-based low-order preconditioning for Poisson problems on rectangular and hexahedral spectral elements has been demonstrated \cite{PedroFischer2019}. Extensions of this approach to finite element problems posed in $\mathds{H}(\mathrm{curl})$ and $\mathds{H}(\mathrm{div})$ spaces, using Nédélec and Raviart–Thomas elements, respectively, have also been reported \cite{Pazner2022deRham}. Related approaches use low-order finite element auxiliary spaces for high-order spectral discretisations \cite{brix2015multilevel}, while weighted multiresolution norm equivalences support robust $p$-version preconditioning \cite{beuchler2004multiresolution}. GPU acceleration and scalability have also been demonstrated for these methods \cite{Phillips2021GPU, Pazner2022GPU}.

Despite its appeal, the LOR approach poses practical challenges in higher-order finite elements implementations. As the polynomial order increases, the resulting discretisation becomes increasingly anisotropic. In particular, \GLL point distributions exhibit grid spacings that scale as $1/P^2$ near element boundaries and $1/P$ in the interior, resulting in strongly non-uniform and anisotropic element geometries at high order \cite{brix2015nested}. This growing anisotropy, together with the resulting irregular element shapes, challenges the core assumptions of standard multigrid and domain-decomposition methods, which rely on moderate, geometry-aligned anisotropy and the effective smoothing of high-frequency error. As a result, the efficiency of these methods can deteriorate significantly for high-order discretisations, particularly on spectrally refined meshes \cite{lottes2005hybrid}. 

At present, there is little support for low-order preconditioning on triangular elements in two dimensions, and no widely available implementations for prismatic, tetrahedral, or pyramidal elements in three dimensions, particularly for industrial-scale simulations on unstructured meshes. Early investigations of triangular finite elements for elliptic problems reported weaker performance for $\mathds{P}1$ discretisations compared to $\mathds{Q}1$ alternatives \cite{deville1994preconditioned}. Low-order preconditioning strategies for simplex elements in two and three dimensions have also been explored using algebraic multigrid within discontinuous Galerkin spectral element solvers \cite{Olson2007}, as well as in the context of overlapping Schwarz methods, where bounded condition number behaviour was observed under both mesh refinement $h$ and polynomial-order refinement $P$ \cite{schoberl2008additive}. More recent work on triangular elements demonstrated improved conditioning through enrichment of the low-order space and the use of alternative node sets \cite{WARBURTON20001, Chalmers2018}. Pazner \cite{pazner2026high} established degree-robust spectral equivalence for a different construction based on Duffy-transformed triangular spaces and a Gauss--Lobatto triangular lattice, including compatibility with quadrilateral elements on mixed meshes. That theory does not directly cover the electrostatic-node triangular and tetrahedral spaces or the prismatic construction considered here. To the authors' knowledge, implementations combining hierarchical high-order bases with LOR preconditioning on prismatic and tetrahedral meshes remain unavailable outside the present work.

In this work, a Low-Order Refined (LOR) preconditioning strategy is developed for the pressure Poisson system arising in the incompressible Navier–Stokes solver of \nekpp. This extends the preliminary formulation in \cite{khurana2024ECCOMAS} through a revised algorithm compatible with conjugate gradient solvers. The proposed approach provides an industrial-scale implementation of LOR that supports mixed two- and three-dimensional element meshes, including simplicial elements, and demonstrates scalability to problem sizes exceeding $10^7$ degrees of freedom. The formulation applies to both Lagrange (\textit{nodal}) and hierarchical (\textit{modal}) basis functions and extends LOR preconditioning to prismatic and tetrahedral meshes. The resulting framework is designed for large-scale computations and leverages algebraic multigrid methods via \nekpp's PETSc interface \cite{petsc-web-page} to enable robust and scalable solution of the associated low-order systems.

%% file: main_methodology.tex
\subsection{Incompressible Navier-Stokes equations}
We consider the LOR preconditioner in the context of solving the incompressible Navier-Stokes equations.
The set of equations describes the motion of an incompressible Newtonian fluid in a Eulerian frame of reference which is governed by the conservation of momentum and mass equations
\begin{align}
    \frac{\partial \mathbf{u}}{\partial t} + \mathbf{u} \cdot \nabla \mathbf{u} + \nabla p - \nu \nabla^2 \mathbf{u} &= 0 \label{eq.navierStokes}, \\
    \nabla \cdot \mathbf{u} &= 0 \label{eq.incompressibility},
\end{align}
where $\mathbf{u}$ is the $d$-dimensional velocity field $(d = 2 \, \text{or} \, 3)$, $p = p / \rho$ the kinematic pressure field (where we assume a unit density $\rho = 1$) and $\nu$ the kinematic viscosity.

We advance the system of equations in time via a velocity correction scheme (VCS) \cite{KaIsOr91}, an operator splitting scheme which combines an explicit-in-time convection, followed by an implicit pressure Poisson problem and $d$ implicit Helmholtz problems \footnote{Helmholtz equations take the form \(-\nabla^{2}u(\mathbf{x})+\lambda u(\mathbf{x})=f(\mathbf{x})\), where \(u\) denotes the solution field and \(\lambda\) is a scalar parameter. For \(\lambda=0\), the equation reduces to the Poisson equation.
} for each component of the velocity vector $V = [u, v, w]^T$.
The pressure problem is elliptic and the pressure serves as an auxiliary scalar to enforce the incompressibility condition in equation \eqref{eq.incompressibility}.
The elliptic nature of the pressure makes the Poisson problem ill-conditioned for large-scale problems and complex, anisotropic computational meshes.

\subsection{The \spectralhp method}
We discretize the computational domain using the \spectralhp element method, a high-order Galerkin method; see \cite{karniadakis2005spectral} for further details.
The \spectralhp method combines the geometric flexibility of the finite element method using mesh elements with size $h$ and the exponential convergence of spectral methods through basis functions with polynomial order $P$.
Accordingly, the computational domain $\Omega$ is partitioned into non-overlapping subdomains $\Omega^e$, referred to as \textit{elements}, which collectively form the computational mesh.
The local subdomains $\Omega^e$ are assembled into the global domain $\Omega$ using an assembly operator $\mathcal{A}^T$ (also referred as scatter/gather) which ensures $C^0$ continuity across subdomains \cite{karniadakis2005spectral}.

In this work, we consider test functions with a hierarchical structure, also referred to as a \textit{modal} basis \cite{VoShKi10}. More specifically, we use the \textit{modified} $C^0$ basis \cite{karniadakis2005spectral}, in which the two linear modes represent the element boundaries and the higher-order modes vanish at the boundaries and are supported in the element interior. We construct the $d$-dimensional basis $\phi(\mathbf{\xi})$ as a tensor product of one-dimensional bases with zeros $\xi$ defined within the standard domain $\mathbf{\xi} \in [-1, 1]^d$ using Gauss or Gauss-Lobatto-Legendre (GLL) point distributions. The integration with Gauss type points, also known as Gaussian quadrature, allows exact evaluation of polynomials up to degree $2Q-1$ using $Q$ zeros and, hence, they are regularly used for higher-order finite element methods.

\subsection{Solving the linear system}
The \spectralhp discretisation of the \IncNS equations using the VCS leads to multiple global linear systems of the form
\begin{equation}
    \mathbf{H}\,\hat{\mathbf{u}} = \mathbf{f},
    \label{eq:helm_system}
\end{equation}
where $\mathbf{H} \in \mathbb{R}^{N_{\mathrm{dof}} \times N_{\mathrm{dof}}}$ 
is the discrete linear operator, 
$\hat{\mathbf{u}} \in \mathbb{R}^{N_{\mathrm{dof}}}$ is the vector of unknown global coefficients, 
and $\mathbf{f} \in \mathbb{R}^{N_{\mathrm{dof}}}$ denotes the discrete forcing term.

The Helmholtz operator is sparse, symmetric and positive-definite and the Poisson operator is sparse, symmetric and indefinite.
We focus our discussion in this work on the pressure Poisson operator.


Within the current \nekpp framework, iterative Krylov subspace methods are available to solve large-scale problems and we consider the Conjugate Gradient \cite{hestenes_methods_1952} and a restarted Generalized Minimal Residual (GMRES) \cite{Saad2003} algorithm.
We measure convergence of the linear solvers with the residual $\hat{\mathbf{r}}_{k} = \mathbf{f} - \mathbf{H}\hat{\mathbf{u}}_{k}$ and compare the residual's $L_2$ norm against a prescribed tolerance $\tau$ as $||\hat{\mathbf{r}}_{k}||_2 < \tau$.
In particular for the elliptic Poisson problem, the ill-conditioning leads to large iteration numbers which become the dominant computational cost of solving the \IncNS equations for complex geometries.

To achieve an equivalent system with better iterative properties, a preconditioning matrix $\mathbf{M}$ is applied such that
\begin{equation}
    \mathbf{M}^{-1}(\mathbf{H}\hat{\mathbf{u}} - \mathbf{f}) = \mathbf{M}^{-1} \hat{\mathbf{r}}  = 0.
    \label{eqn:precon_nektar}
\end{equation}
In practice, the preconditioner is applied to the residual at every iteration to compute the preconditioned direction $\hat{\mathbf{z}}_{k} = \mathbf{M}^{-1} \hat{\mathbf{r}}_{k}$. For $\mathbf{M}$ to be effective, the condition number $\kappa$ of the modified system must be significantly reduced, such that $\kappa(\mathbf{M}^{-1}\mathbf{H}) \ll \kappa(\mathbf{H})$ and ideally approaches unity.

\subsection{LOR preconditioner}
\begin{figure}[H]
    \centering
    \includegraphics[width=0.5\textwidth]{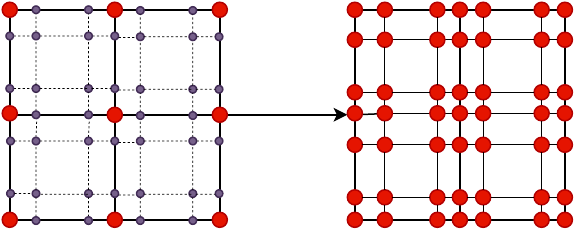}
    \caption{Conversion of four higher-order nodal finite elements at \textit{P}=3 to the corresponding \textit{P}=1 low-order refined (LOR) representation. The LOR grid is constructed by creating \textit{P}=1 elements at the locations of the quadrature points of the higher-order element. Red dots denote the linear mesh nodes, while purple dots indicate the quadrature points within each higher-order element. The resulting LOR mesh contains only linear elements and preserves the GLL-based point distribution of the higher-order discretisation.}
    \label{fig:LOR_concept}
\end{figure}

For the formulation of the LOR preconditioner matrix $\mathbf{M}_{LOR}$, a separate low-order (\textit{P}=1) ``refined'' finite element system $\mathbf{H}_{L}$ is constructed that is \textit{spectrally equivalent} to the higher-order \spectralhp system $\mathbf{H}$. Previous LOR formulations have primarily considered rectangular and hexahedral meshes with tensor-product Lagrange polynomial bases defined on GLL quadrature points \cite{PedroFischer2019, Pazner2022deRham}. In the present implementation, the higher-order discretisation instead uses the modified $C^0$ modal basis described previously, which is interpolated onto the nodal LOR space. The conceptual representation of this conversion is shown in Fig.~\ref{fig:LOR_concept}, where the LOR mesh on the right is a sub-structured g

For the high-order linear system $\mathbf{H}$, every internal DOF is coupled with every other DOF in an element. The conversion from higher-order modal to the LOR space $\mathbf{H}_{L}$ increases the sparsity as now the DOF are only connected to their immediate neighbours. The difference in connectivity between the two discretisations is further accentuated in mathematical operators, even if the number of DOF is the same (as seen in Fig.~\ref{fig:LOR_concept}). Thus, the operators on the lower-order discretisation are cheaper to evaluate than their higher-order counterparts. 

The construction of the LOR preconditioner proceeds in two main steps: first, the interpolation of the high-order element space to its low-order refined equivalent; and second, the definition of the LOR finite element space for the different element types considered. Each step is described in detail in the following sections.

\subsubsection{Interpolation to LOR space}
\label{subsec:lor_interpolation}
The interpolation is presented as two transformation steps: a) a \textit{Vandermonde} transformation that moves from a modal to a nodal expansion basis in the higher-order space, and b) a collocation projection from the higher-order nodal space to the $P=1$ LOR space. 

To interpolate any basis function to the LOR space for the LOR preconditioner, we distinguish between the three bases:
\begin{itemize}
    \item \textbf{HO modal}: the modified $C^{0}$ hierarchical basis
    \cite{karniadakis2005spectral}, $\phi_p(\boldsymbol{\xi}) \in \mathcal{P}^{P}$,
    whose two boundary modes are linear and whose interior modes are Jacobi
    polynomials $J^{1,1}_{p-1}$ weighted by
    $\tfrac{1-\xi}{2}\tfrac{1+\xi}{2}$, with coefficients
    $\hat{\mathbf{u}} = \{\hat{u}_p\}_{p=0}^{P}$,
    \item \textbf{HO nodal}: Lagrange basis $h_p(\boldsymbol{\xi}) \in \mathcal{P}^{P}$
    defined at the $P+1$ Gauss-Lobatto-Legendre (GLL) nodes
    $\boldsymbol{\xi} = \{\xi_j\}_{j=0}^{P}$, with coefficients
    $\bar{\mathbf{u}} = \{\bar{u}_p\}_{p=0}^{P}$,
    \item \textbf{LOR}: continuous piecewise-linear basis
    $l_p(\boldsymbol{\xi})$, locally in $\mathcal{P}^{1}$ on each of the $P$
    sub-intervals of the refined grid defined by $\boldsymbol{\xi}$, with
    coefficients $\hat{\mathbf{u}}' = \{\hat{u}'_p\}_{p=0}^{P}$.
\end{itemize}

\begin{figure}[H]
\centering
\begin{subfigure}{0.32\textwidth}
  \centering
  \includegraphics[width=\linewidth]{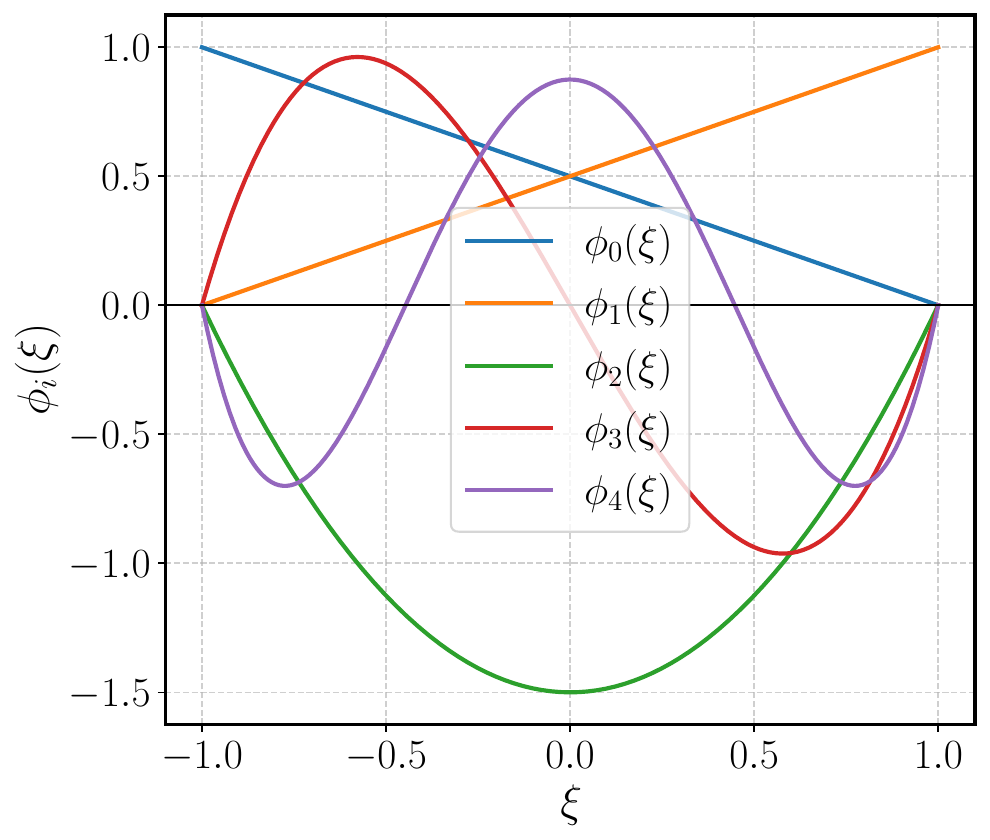}
  \caption{Modified Jacobi basis}
  \label{subfig:modal_P4}
\end{subfigure}
\begin{subfigure}{0.32\textwidth}
  \centering
  \includegraphics[width=\linewidth]{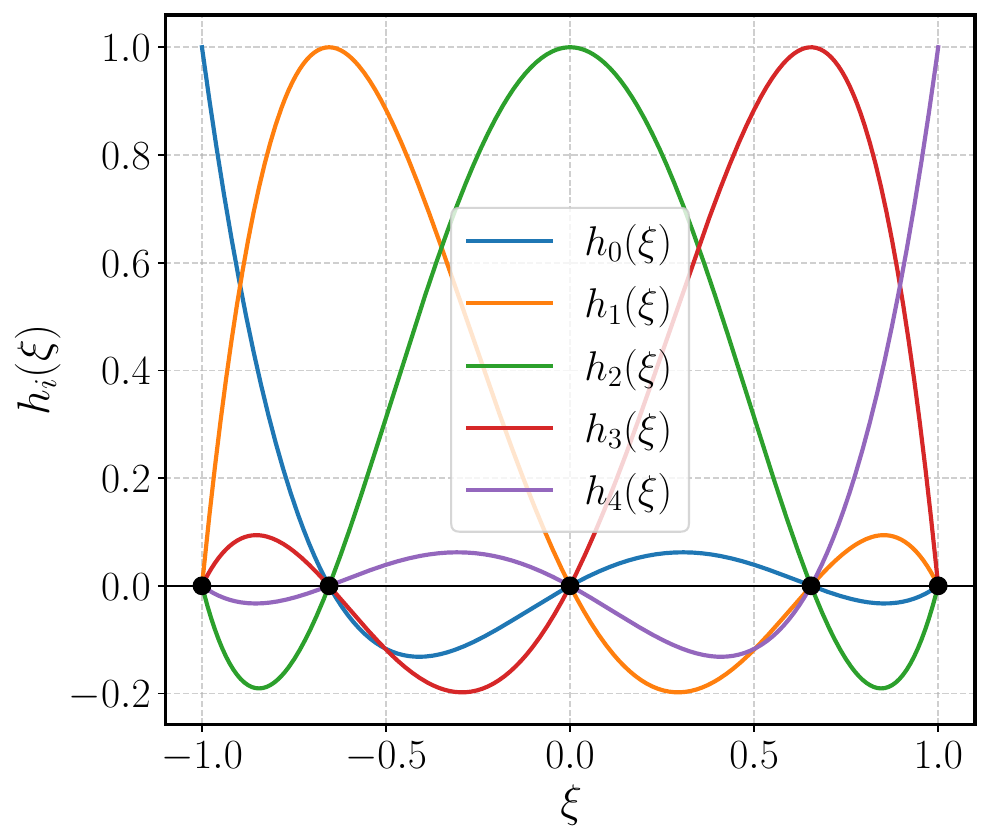}
  \caption{Nodal Lagrange basis}
  \label{subfig:nodal_P4}
\end{subfigure}
\begin{subfigure}{.32\textwidth}
  \centering
  \includegraphics[width=\linewidth]{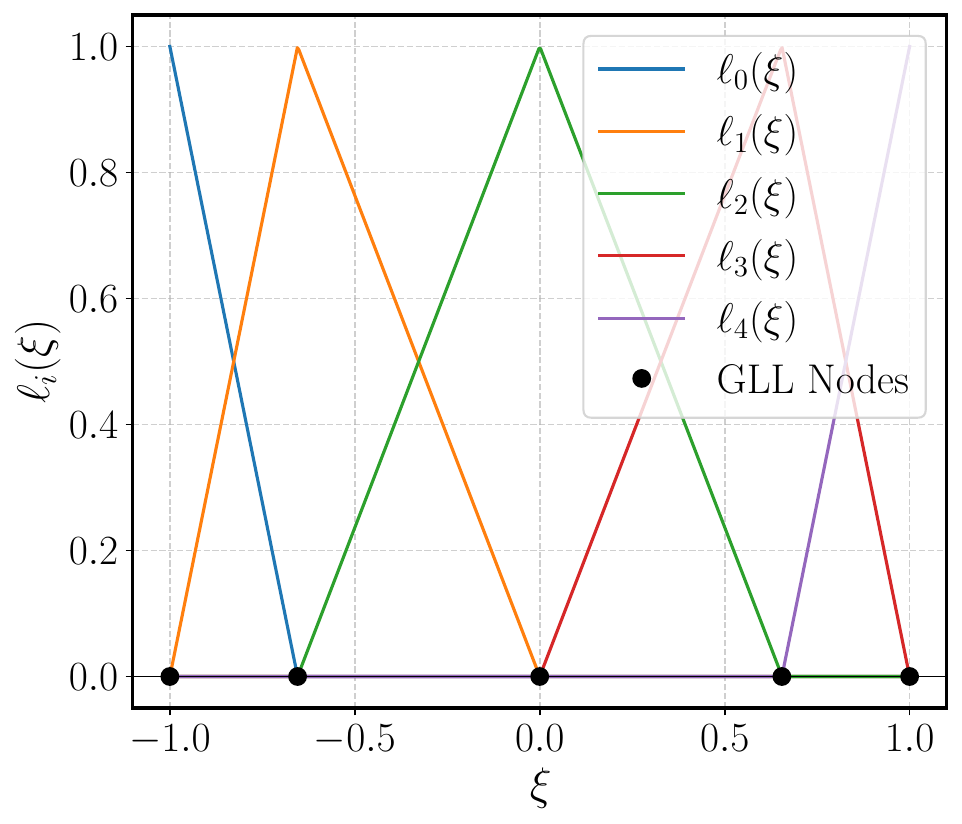}
  \caption{LOR piecewise basis}
  \label{subfig:LOR_nodal_P4}
\end{subfigure}
\caption{
Comparison of the 1D \nekpp basis functions and the corresponding low-order refined (LOR) basis for \(P=4\). In the LOR construction, the original high-order element is divided into \(n_{\mathrm{split}}=4\) sub-intervals, and each basis function is represented as a continuous piecewise-linear function over this refined partition. This figure illustrates how the HO function space is represented in the LOR setting via a refined nodal discretisation using the same \(P+1\) nodal points.}
\label{fig:LOR_vs_nodal_vs_modal_P4}
\end{figure}

Fig.~\ref{fig:LOR_vs_nodal_vs_modal_P4} shows the three bases at $P=4$: the
modal and nodal high-order bases, and the piecewise-linear LOR basis defined
on the sub-grid of the same GLL nodes. The corresponding polynomial representations of the solution in the higher order space are
\begin{equation}
    \underbrace{\sum_{p=0}^{P} h_p(\boldsymbol{\xi})\, \bar{u}_p}_{\text{nodal}} = u^{P}(\boldsymbol{\xi}) = \underbrace{\sum_{p=0}^{P} \phi_p(\boldsymbol{\xi})\, \hat{u}_p}_{\text{modal}}
    \label{eq:uLagrange}
\end{equation}

Each subinterval of the original high-order element is represented by a local linear Lagrange basis \(l_p^{(e)}(\boldsymbol{\xi})\), such that \(l_p^{(e)} \in\mathcal{P}^1\) within each sub-element $\Omega^e$. We denote the number of subintervals used to partition the high-order element by \(n_{\mathrm{split}}\); for a one-dimensional element with (P+1) nodal points, \(n_{\mathrm{split}}=P\).

The LOR expansion over the original high-order element is then assembled as the continuous piecewise-linear representation
\begin{equation}
  u^{P=1}(\boldsymbol{\xi}) \;=\; \sum_{p=0}^{P} l_p(\boldsymbol{\xi})\,\hat{u}'_p,
  \label{eq:uLOR}
\end{equation}
while the index range reflects the \(P+1\) nodal degrees of freedom inherited from the HO nodal set. Each basis \(l_p\) is therefore piecewise linear, being locally in \(\mathcal{P}\) for each sub-interval $n_{split}$, as visualised in Fig.~\ref{subfig:LOR_nodal_P4}.

A polynomial can be represented by any other basis of the same polynomial order in the same space using the \textit{Vandermonde transformation} \cite{WARBURTON20001}. The standard Vandermonde matrix $\mathcal{V} \in\mathbb{R}^{(P+1)\times (P+1)}$ relates the HO modal coefficients $\hat{\mathbf{u}}$ and the HO nodal coefficients $\bar{\mathbf{u}}$ as
\begin{equation}
    \bar{\mathbf{u}} = \mathcal{V} \hat{\mathbf{u}}.
    \label{eqn:modaltonodal_coeff}
\end{equation}
Here, $\mathcal{V}_{ij} = \phi_j(\boldsymbol{\xi}_i)$, where the invertibility of $\mathcal{V}$ is guaranteed for a valid set of basis functions and distinct nodes. Throughout, $\Phi, \mathbf{h}, \boldsymbol{\ell} \in
\mathbb{R}^{Q \times (P+1)}$ denote the modal, nodal and LOR basis functions
evaluated at the $Q$ quadrature points. Hence,
$u^{P} = \mathbf{h}\bar{\mathbf{u}} = \Phi\hat{\mathbf{u}}$ at these points.
We take $Q = P+1$, with the quadrature points collocated with the GLL nodes.
Using Eqn.~\eqref{eqn:modaltonodal_coeff} and requiring this equality for arbitrary $\hat{\mathbf{u}}$ gives
\begin{equation}
\mathbf{h}\mathcal{V} = \Phi
\quad \implies \quad
\mathbf{h} = \Phi\mathcal{V}^{-1}.
\label{eqn:modaltonodal_basis}
\end{equation}

To map the smooth high-order nodal solution $u^P(\boldsymbol{\xi})$ (Eqn.~\eqref{eq:uLagrange}) onto the low-order representation $u^{P=1}(\boldsymbol{\xi})$ (Eqn.~\eqref{eq:uLOR}), we employ a \textit{collocation projection}. This imposes the constraint that both expansions coincide exactly at the set of $P+1$ quadrature nodes $\boldsymbol{\xi}$. Since both the HO Lagrange basis $\{h_p\}$ and the piecewise-linear LOR basis $\{l_p\}$ satisfy the Kronecker delta property at these nodes ($h_p(\xi_j) = l_p(\xi_j) = \delta_{pj}$), the expansions are uniquely determined by their values at $\boldsymbol{\xi}$. Consequently, the expansion coefficients in both spaces are identical, such that $\hat{\mathbf{u}}' \equiv \bar{\mathbf{u}}$. The Kronecker delta property likewise gives $\boldsymbol{\ell} = \mathbf{h}$ at the collocated points, so the nodal and LOR representations differ only
in the polynomial space each basis spans.

The transformation between the modal coefficients $\hat{\mathbf{u}}$ and the nodal/LOR coefficients $\mathbf{u}'$ is given by the Vandermonde matrix $\mathcal{V}$:
\begin{equation}
\hat{\mathbf{u}}' = \mathcal{V} \hat{\mathbf{u}}.
\label{eqn:modaltoLOR_coeff}
\end{equation}

\subsubsection{Operator Equivalence and Preconditioning} To formulate the preconditioner, we must map the linear system operators between these spaces. From Eqn. \eqref{eq:helm_system}, the discrete HO system is $\mathbf{H}\hat{\mathbf{u}} = \hat{\mathbf{r}}$, where $\mathbf{H}$ is the modal Helmholtz matrix and the residual vector $\hat{\mathbf{r}}$ represents the discrete inner products $(\phi_p, f)$ of the forcing function $f$ with the modal basis. The residual vectors in the two bases are obtained by taking the inner product of $f$ with the corresponding basis functions:
$\hat{\mathbf{r}} = \Phi^{\mathsf{T}}\mathbf{W}f$ and
$\bar{\mathbf{r}} = \mathbf{h}^{\mathsf{T}}\mathbf{W}f$, where
$\mathbf{W}$ contains the quadrature weights scaled by the Jacobian.
Using $\mathbf{h} = \Phi\mathcal{V}^{-1}$ gives
\begin{equation}
\bar{\mathbf{r}} = (\Phi\mathcal{V}^{-1})^{\mathsf{T}}\mathbf{W}f
= \mathcal{V}^{-\mathsf{T}}\hat{\mathbf{r}}.
\label{eqn:residual_transform}
\end{equation}
Thus, coefficients transform with $\mathcal{V}$, whereas residuals transform with
$\mathcal{V}^{-\mathsf{T}}$. This preserves the duality pairing,
\[
\bar{\mathbf{u}}^{\mathsf{T}}\bar{\mathbf{r}} =
\hat{\mathbf{u}}^{\mathsf{T}}\hat{\mathbf{r}}.
\]

Eq.~\eqref{eqn:residual_transform} demonstrates that only a transformation of the coefficient-space residuals $\hat{\mathbf{r}}$ is needed to interpolate from the HO modal space to the LOR space. The matrix $\mathcal{V}^{-\mathsf{T}}$ is calculated once and applied to $\hat{\mathbf{r}}$ at each iteration. Furthermore, this is an invertible operation that ensures a consistent mapping to and from the HO and LOR spaces.


As established in Sec.~\ref{subsec:lor_interpolation}, the collocation projection identifies the HO nodal and LOR coefficient vectors at the common GLL nodes. Hence, for the collocation transfer used here, $\hat{\mathbf{u}}'=\bar{\mathbf{u}}$ and $\hat{\mathbf{r}}'=\bar{\mathbf{r}}$, so that the LOR system $\mathbf{H}_L\hat{\mathbf{u}}'=\hat{\mathbf{r}}'$ can be written in the corresponding nodal coordinates as $\mathbf{H}_L\bar{\mathbf{u}}=\bar{\mathbf{r}}$. Substituting the coefficient transformation \eqref{eqn:modaltonodal_coeff} and residual transformation \eqref{eqn:residual_transform} into this LOR system yields:

\begin{align}
    \mathbf{H}_L (\mathcal{V} \hat{\mathbf{u}}) &= \mathcal{V}^{-T} \hat{\mathbf{r}} \nonumber \\
    \implies \quad \mathcal{V}^{T} \mathbf{H}_L \mathcal{V} \hat{\mathbf{u}} &= \hat{\mathbf{r}}.
\end{align}
Comparing this result to Eq.~\eqref{eq:helm_system}, we identify the LOR preconditioner projected into the modal space, $\mathbf{M}_{LOR}$, as:
\begin{equation}
    \mathbf{M}_{LOR} = \mathcal{V}^{T}\,\mathbf{H_L}\,\mathcal{V}.
    \label{eqn:M_LOR_modal}
\end{equation}
For implementations using a HO \textit{nodal} basis (where $\mathcal{V}=\mathbf{I}$), this simplifies to
\begin{equation}
    \mathbf{M}_{LOR} = \mathbf{H_L}.
    \label{eqn:M_LOR_nodal}
\end{equation}
The two definitions in Eqns.~\eqref{eqn:M_LOR_modal}
and~\eqref{eqn:M_LOR_nodal} describe the same preconditioned system in
different bases. Coefficients transform as
$\bar{\mathbf{u}} = \mathcal{V}\hat{\mathbf{u}}$ and residuals as
$\bar{\mathbf{r}} = \mathcal{V}^{-\mathsf{T}}\hat{\mathbf{r}}$. Applying
both to the nodal system $\mathbf{H}_{\mathrm{nodal}}\bar{\mathbf{u}} =
\bar{\mathbf{r}}$ recovers the modal system
$\mathbf{H}_{\mathrm{modal}}\hat{\mathbf{u}} = \hat{\mathbf{r}}$ with
\begin{equation}
    \mathbf{H}_{\mathrm{modal}}
    = \mathcal{V}^{\mathsf{T}}\mathbf{H}_{\mathrm{nodal}}\mathcal{V},
    \label{eqn:H_congruence}
\end{equation}
which is the same congruence that relates $\mathbf{M}_{LOR}$ to
$\mathbf{H}_{L}$ in Eqn.~\eqref{eqn:M_LOR_modal}. The two factors of
$\mathcal{V}$ then cancel in the preconditioned operator,
\begin{equation}
    \mathbf{M}_{LOR}^{-1}\mathbf{H}_{\mathrm{modal}}
    = \left(\mathcal{V}^{\mathsf{T}}\mathbf{H}_{L}\mathcal{V}\right)^{-1}
      \left(\mathcal{V}^{\mathsf{T}}\mathbf{H}_{\mathrm{nodal}}
      \mathcal{V}\right)
    = \mathcal{V}^{-1}\left(\mathbf{H}_{L}^{-1}
      \mathbf{H}_{\mathrm{nodal}}\right)\mathcal{V}.
    \label{eqn:LOR_similarity}
\end{equation}
The modal and nodal preconditioned operators are related by a similarity
transformation and therefore have identical spectra. Krylov convergence is
consequently independent of the choice of high-order basis, and any
difference observed in practice is attributable to the conditioning of
$\mathcal{V}$ itself.

%% file: algorithm.tex
The LOR preconditioner algorithm comprises two distinct stages: a one-time \textit{Build stage} and a per-iteration \textit{Do stage}. The \textit{Build stage} handles the initial construction of the LOR mesh and its associated matrix systems, which are then reused by the \textit{Do stage} during the iterative solve. The setup costs during the \textit{Build stage} can be substantial, especially for large, complex 3D meshes, but they are amortised by the gains in iterative performance of the \textit{Do stage} due to improved conditioning properties. This section details the implementation of each stage, beginning with a high-level overview of the complete algorithm.

\paragraph{Assembly maps and multiplicity}
The LOR preconditioner repeatedly transfers vectors between global and
element-local representations. We define the assembly operators as follows: the scatter operator $\mathcal{A}$ maps global degrees of freedom to element-local degrees of freedom, while the gather operator $\mathcal{A}^T$ maps element-local contributions back to the global space. Since neighbouring elements share interface DOFs, the gather operation introduces repeated contributions at the global level. To account for the resulting over-counting, a diagonal multiplicity matrix $\mathbf{D}$ is defined such that the consistent global recovery is given by
\begin{equation}
    u_g = \mathbf{D}^{-1}\mathcal{A}^T u_l.
\end{equation}
Here, $\mathbf{D} = \mathrm{diag}(d_1, d_2, \dots, d_{N_{\mathrm{dof}}})$, where each $d_i$ denotes the number of local DOFs associated with global DOF $i$. In the following, subscripts $H$ and $L$ distinguish the corresponding
assembly maps and multiplicity matrices in the high-order and LOR spaces.

\paragraph{Preconditioner application} Fig.~\ref{fig:LOR_process} provides an overview of one application of the LOR
preconditioner. The input is the high-order (HO) residual vector
$\hat{\mathbf{r}}_{H,g}$ in the global coefficient space. The residual is first
scattered to the element-local HO space using $\mathcal{A}_H$ and then
interpolated to the corresponding LOR space. This interpolation, detailed in
Section~\ref{subsec:lor_interpolation}, depends on the chosen HO expansion
basis and may require construction of $\mathcal{V}$.

The resulting LOR vector is transferred to the representation required by the
inner solver and the LOR system is solved using the chosen
\textit{inner-iteration strategy}. The resulting correction is then projected
back to the HO element space and gathered into the global HO space. The
assembly operations in the HO and LOR spaces are described by
$\mathcal{A}_H$ and $\mathcal{A}_L^e$, respectively, with the corresponding
multiplicity matrices $\mathbf{D}_H$ and $\mathbf{D}_L^e$ accounting for
shared DOFs.

\input{algos/LOR_process_global}

The \textit{Do stage} is detailed in Algorithm~\ref{alg:DoLORPreconditioner}. Its central operation is the globally assembled LOR solve in Step~5, which maps $\hat{\mathbf{r}}_{L,g}$ to $\hat{\mathbf{z}}_{L,g}$. The resulting LOR correction is transformed to the HO element space in Step~6 and assembled into the global HO correction $\hat{\mathbf{z}}_{H,g}$ in Step~7, which is returned to the outer solver (CG/GMRES).

\input{algos/doLOR}



The two multiplicity matrices account for shared degrees of freedom at different stages of the transfer. In the HO space, $\mathbf{D}_H=\mathcal{A}_H^T\mathcal{A}_H$ counts the number of element-local copies associated with each global HO degree of freedom. Step~1 divides each global residual entry by this multiplicity before $\mathcal{A}_H$ scatters the residual to the HO element-local space.

Step~3 applies the inverse transpose of the generalized Vandermonde
operator, $\mathcal{V}^{-T}$, to express the HO element-local residual in
the LOR basis on the refined mesh, yielding the assembled global vector $\hat{\mathbf{r}}_{L,g}$. After the transformation in Step~3, $\hat{\mathbf{r}}_{L,g}$ lies in the assembled global LOR coefficient space represented by the rightmost mesh in Fig.~\ref{fig:LOR_process}. During the Build stage, the elemental LOR assembly map $\mathcal{A}_L^e$ defines the corresponding multiplicity matrix
\begin{equation}
    \mathbf{D}_L^e
    =
    (\mathcal{A}_L^e)^T\mathcal{A}_L^e.
\end{equation}
Each diagonal entry of $\mathbf{D}_L^e$ counts the refined subelements that share the corresponding assembled LOR degree of freedom. Step~4 applies
$(\mathbf{D}_L^e)^{-1}$ to account for these repeated subelement contributions.

The correction $\mathbf{D}_L^e$ is essential for simplicial and mixed-element meshes because the mapping from HO local coefficients to LOR degrees of freedom is not one-to-one. Multiple LOR degrees of freedom can reference the same HO-local coefficient, which causes systematic over-counting in the LOR residual unless element-wise multiplicity is applied. For tensor-product elements, such as quadrilaterals and hexahedra, this correction is typically trivial. For simplicial elements, it is critical. The LOR system is then solved directly in the global coefficient space in Step~5. 

Equation~\ref{eqn:LOR_alloperations} expresses the LOR preconditioner as a composite linear operator acting on the high-order residual, consisting of inverse multiplicity scaling, assembly, and projection between high-order and low-order refined spaces, and the application of the low-order inverse operator.

\begin{equation}
\hat{\mathbf z}_{H,g}
=
\underbrace{\mathcal{A}_H^T}_{\substack{\text{HO}\\\text{assemble}}}
\underbrace{\mathcal{V}^T}_{\substack{\text{transform}\\\text{to HO}}}
\underbrace{\mathbf M_{\mathrm{LOR}}^{-1}}_
            {\substack{\text{global}\\\text{LOR solve}}}
\underbrace{(\mathbf D_L^e)^{-1}}_
            {\substack{\text{LOR}\\\text{scale}}}
\underbrace{\mathcal{V}^{-T}}_
            {\substack{\text{transform}\\\text{to LOR}}}
\underbrace{\mathcal{A}_H}_
            {\substack{\text{HO}\\\text{scatter}}}
\underbrace{\mathbf D_H^{-1}}_
            {\substack{\text{HO}\\\text{scale}}}
\hat{\mathbf r}_{H,g}.
\label{eqn:LOR_alloperations}
\end{equation}

\paragraph{One-time build stage} The needed operators and mappings in \textit{Do stage} i.e. $\mathcal{A}_H$, $\mathcal{A}_L^e$ and $\mathbf{M_{LOR}}$ are constructed and cached during \textit{Build stage}. The \textit{Build stage} algorithm is detailed in Algorithm \ref{alg:BuildLORPreconditioner}. This stage is characterized by extensive reuse of high-order assembly maps and mesh graph information.

A new mesh graph for the LOR space $\mathcal{M}_L$ is generated (Step 2) by subdividing the original HO mesh graph $\mathcal{M}_H$. A degree-$P$ expansion has $P+1$ nodal points along each tensor-product edge, and hence $n_{\mathrm{split}}=P$ refined subintervals. The subdivision points follow the nodal distribution of the HO element. Two distinct mesh graphs are maintained for the HO and LOR discretisations. The construction for the supported element types is summarised in Section~\ref{subsec:lor_space_construct}, with implementation-level algorithms given in Appendix~\ref{app:lor_mesh_construction}.

\input{algos/buildLOR}

Next, the existing HO assembly map ($\mathcal{A}_H$) is leveraged to generate the element-wise LOR assembly map ($\mathcal{A}_L^e$) and the coefficient mapping $\mathcal{P}_{H \rightarrow L}$ (Step 3). The mapping $\mathcal{P}_{H \rightarrow L}$ maps the indices of the subdivided LOR elements to the parent HO element, along with vertex offsets to be used later for data transfer. After assembling the elemental LOR matrix $\mathbf{H}_L$ using the point locations obtained in $\mathcal{M}_L$ (Step 4), the stage concludes by computing the final preconditioner matrix $\mathbf{M}_{\text{LOR}}$. For the case of an HO modal expansion basis, the generalised Vandermonde matrix $\mathcal{V}$ is also needed. The \textit{Build stage} for the LOR preconditioner is much more time-consuming than any other \nekpp preconditioner. However, it provides a bridge between the HO and LOR spaces, completing the setup required for rapid application in the subsequent \textit{Do stage}.

\subsection{Choice of LOR space}
\label{subsec:lor_space_construct}
The definition of the LOR space is a key component of the preconditioner, because the nodal discretisation used in that space directly affects the conditioning of the resulting auxiliary system. For tensor-product elements, this construction follows naturally from collocated \GLL nodes and a Lagrange basis. Extending the construction to simplicial and prismatic elements is important because these elements provide geometric flexibility for body-fitted meshes of complex and industry-relevant geometries \cite{NekMeshGREEN2024}. For simplicial elements, the absence of a tensor-product nodal structure requires an alternative nodal basis spanning the same polynomial space. In this work, the simplicial LOR space uses electrostatic points $\boldsymbol{\xi}_{\mathrm{ele}}$ in the element interior \cite{HestavenElectro1998}, while boundary points are constrained to the distributions used on neighbouring faces and edges. This produces compatible nodal representations for low-order refinement of triangles and tetrahedra; the prismatic construction combines the triangular subdivision with a one-dimensional GLL subdivision. Appendix~\ref{app:lor_mesh_construction} gives the corresponding connectivity algorithms for two-dimensional elements, tetrahedra, and prisms.

\subsubsection{LOR space for 2D element types}
\begin{figure}[H]
    \centering
    \includegraphics[width=0.5\textwidth]{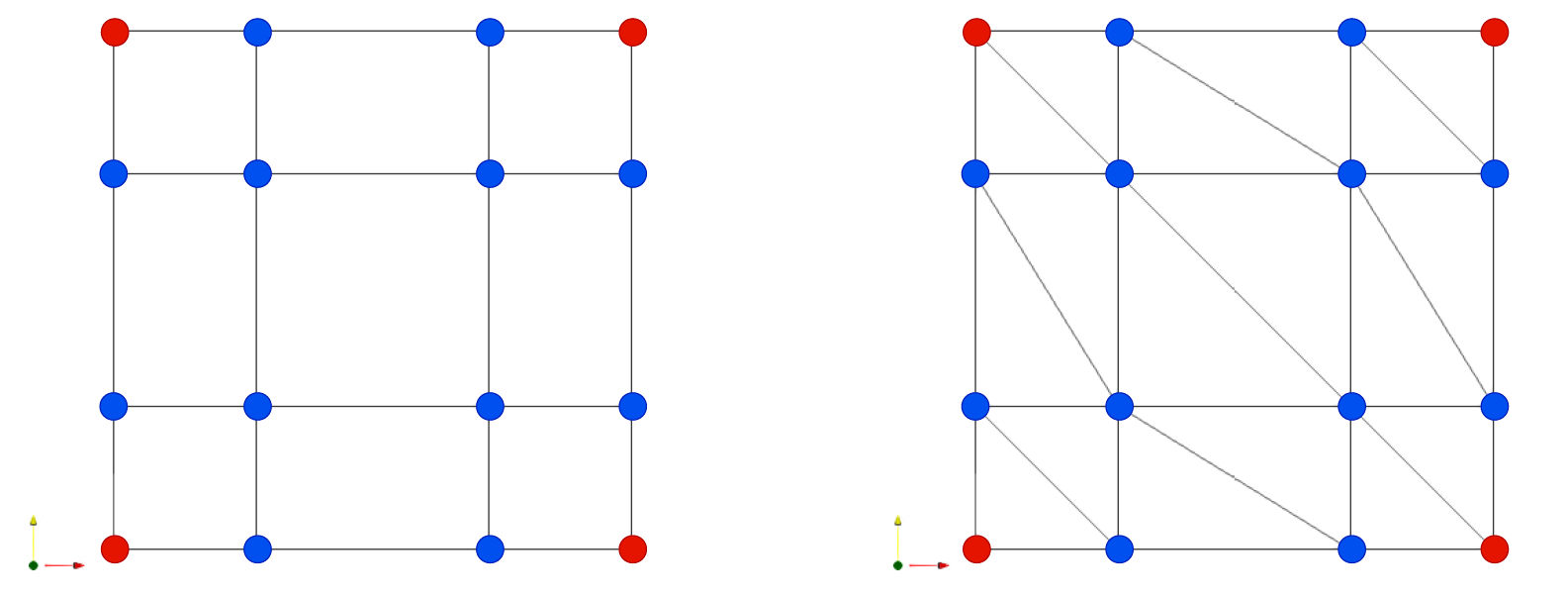} 
    \caption{Available options for LOR refinement of a higher-order quadrilateral element in 2D. The red dots correspond to the vertices of the original quadrilateral element with \textit{P}=3, and the blue dots are interior points in GLL distribution which correspond to vertices in the LOR space. Left: $\mathds{Q}1$ configuration, Right: $\mathds{P}1$ configuration.}
    \label{fig:P1vsQ1}
\end{figure}

There are two options to discretise the quadrilaterals in the LOR space following the work of Canuto \cite{Canuto2009}, namely $\mathds{Q}1$ and $\mathds{P}1$ elements. The differences can be seen in Fig.~\ref{fig:P1vsQ1}, where for a given quadrilateral element at a fixed polynomial order the GLL points in the LOR space can be formulated with $\mathds{Q}1$ (quadrilateral) elements or $\mathds{P}1$ (triangular) elements. 

\begin{figure}[H]
    \centering
    \includegraphics[width=0.5\textwidth]{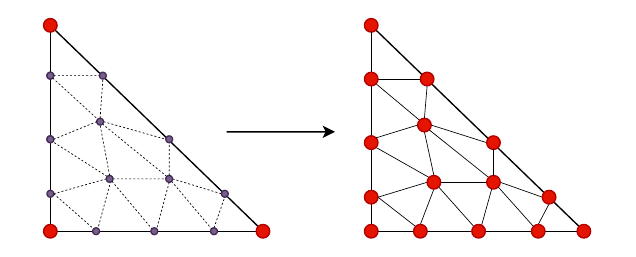} 
    \caption{Conceptual representation of the LOR space for triangles. The state on the left is the higher order triangular element with nodal expansion basis collocated to the electrostatic points $\mathbf{\xi_{ele}}$ at \textit{P}=4. On the right, the LOR space is constructed as a collocated Lagrange basis at \textit{P}=1 on $\mathbf{\xi_{ele}}$ points indicated by red dots.}
    \label{fig:LOR_Triangles}
\end{figure}

Notice that there can be many possible orientations for $\mathds{P}1$ elements, and the choice in Fig.~\ref{fig:P1vsQ1} is one such possibility. For the  $\mathds{Q}1$ configuration,  $n_{\text{split}}^2$ subdivided quadrilaterals are obtained, while the  $\mathds{P}1$ configuration results in $2n_{\text{split}}^2$ subdivided triangles. There are no such options for triangular elements in the LOR space as they can only be discretised into further smaller $\mathds{P}1$ elements. 

\subsubsection{LOR space for 3D element types}
\begin{figure}[H]
\centering
\begin{subfigure}{0.5\textwidth}
  \centering
  \includegraphics[width=\linewidth]{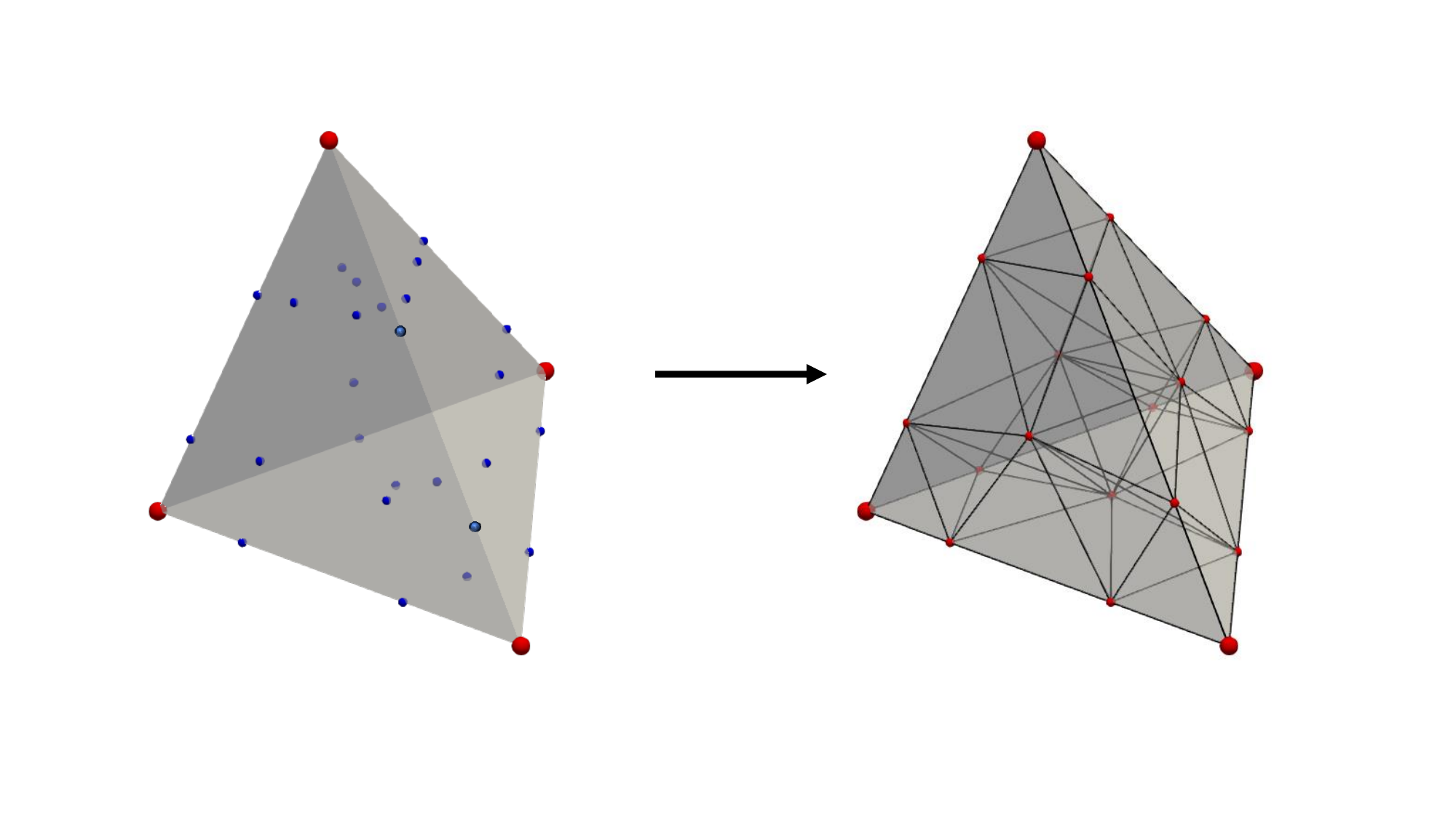}
  \caption{Tetrahedron at \textit{P} = 3}
  \label{subfig:Tet_P3}
\end{subfigure}%
\begin{subfigure}{.5\textwidth}
  \centering
  \includegraphics[width=\linewidth]{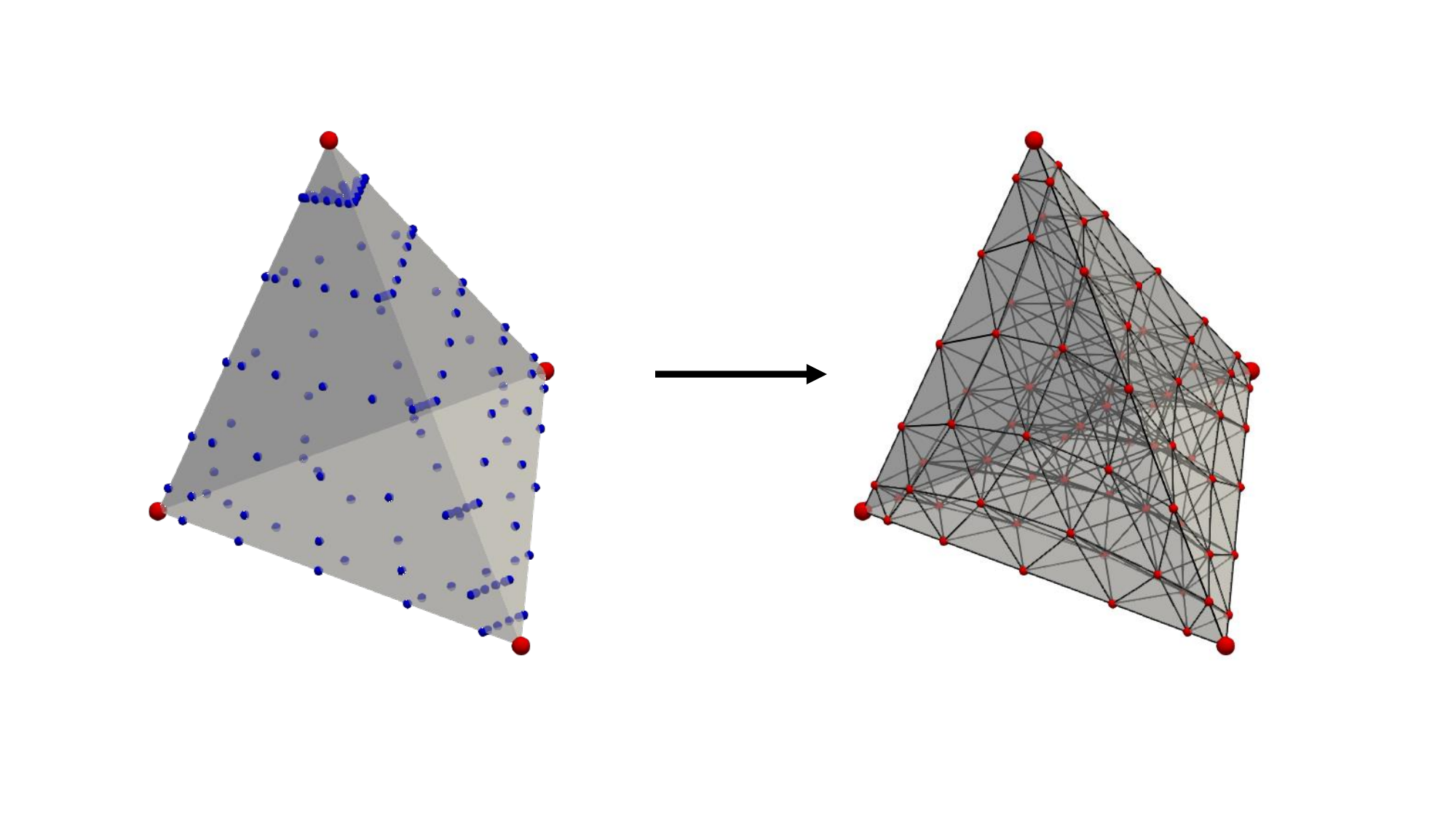}
  \caption{Tetrahedron at \textit{P} = 6}
  \label{subfig:Tet_P6}
\end{subfigure}

\begin{subfigure}{0.5\textwidth}
  \centering
  \includegraphics[width=\linewidth]{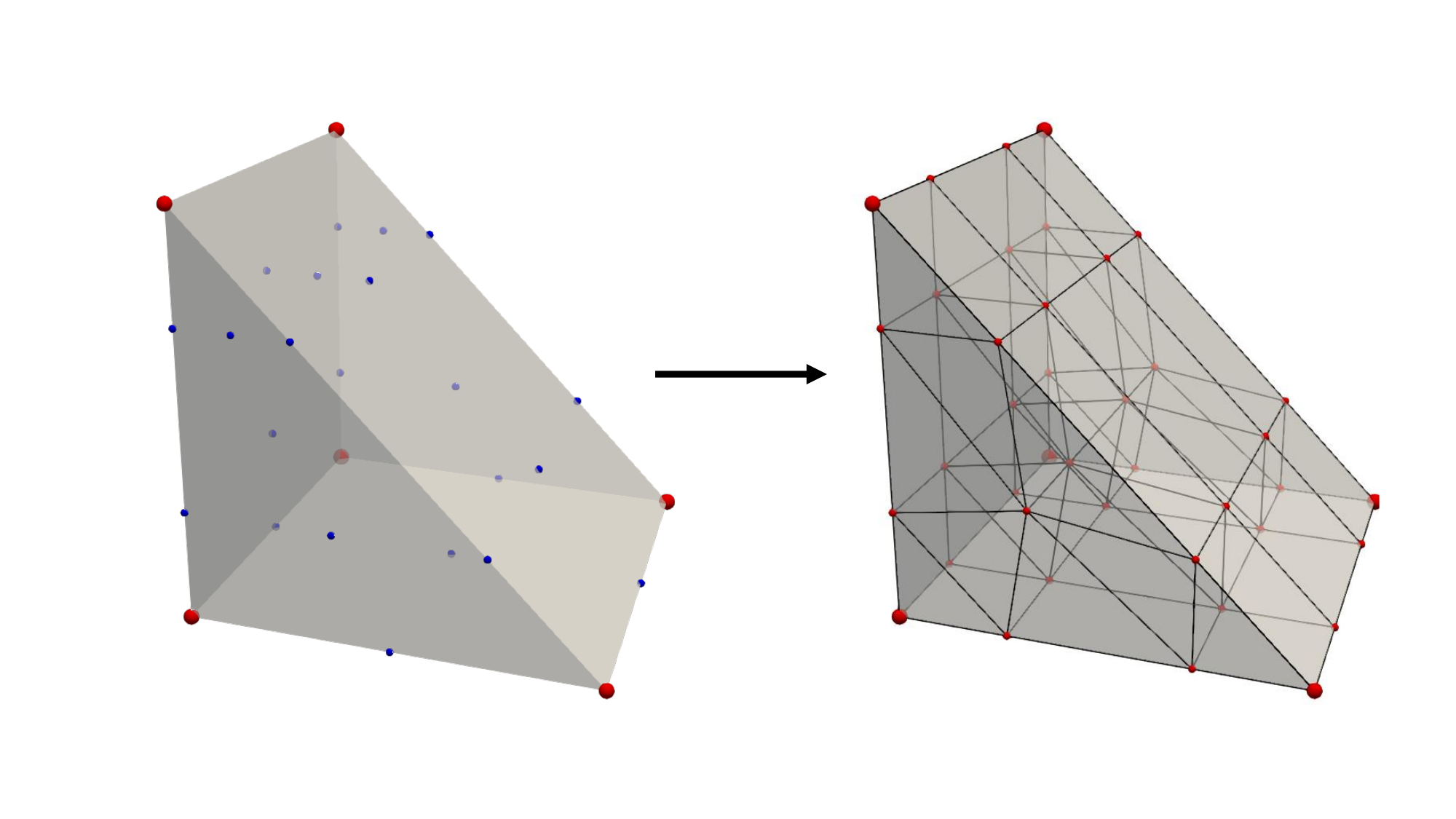}
  \caption{Prism at \textit{P} = 3}
  \label{subfig:Prism_P3}
\end{subfigure}%
\begin{subfigure}{.5\textwidth}
  \centering
  \includegraphics[width=\linewidth]{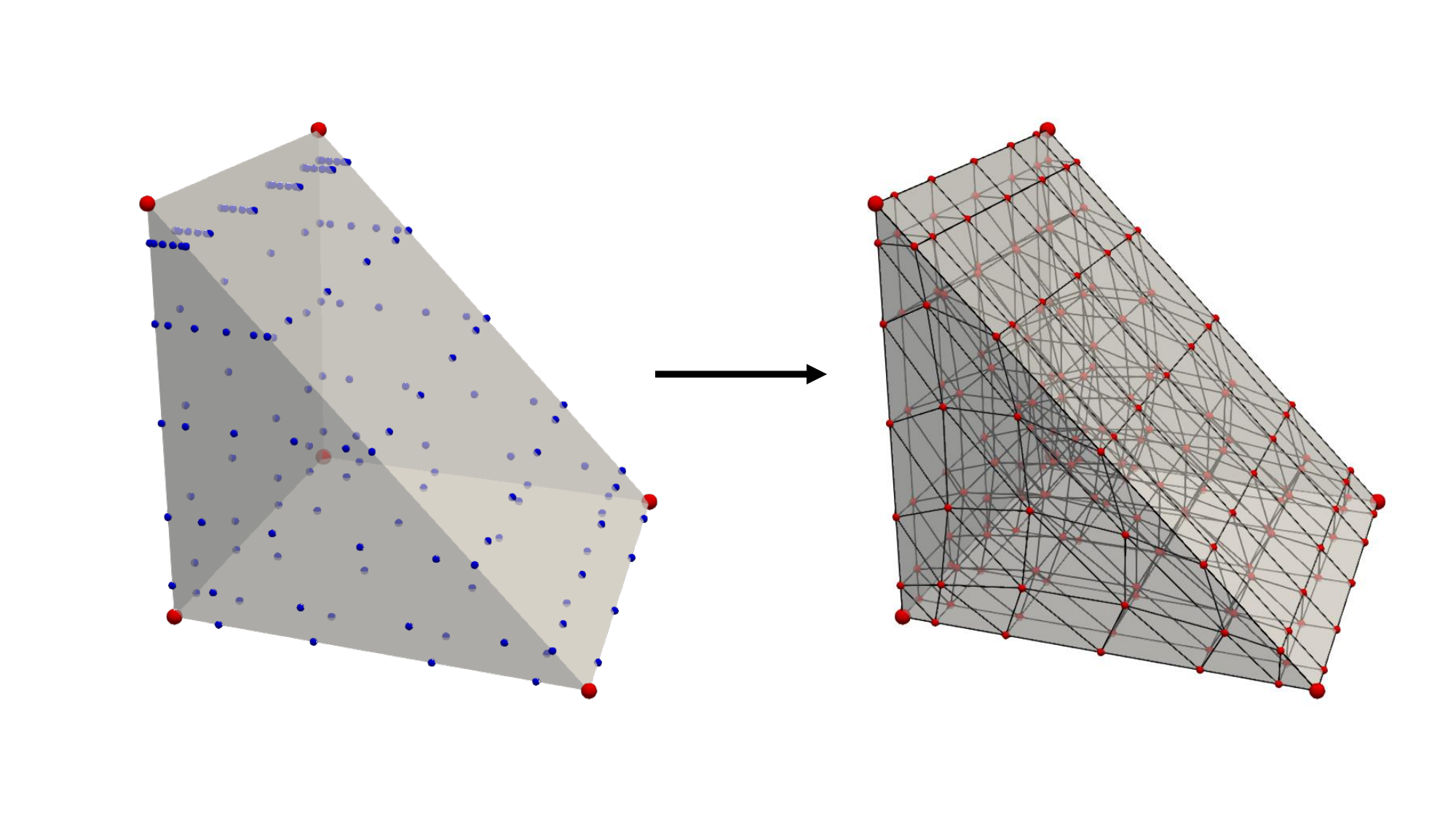}
  \caption{Prism at \textit{P} = 6}
  \label{subfig:Prism_P6}
\end{subfigure}
\caption{Visual representation of the LOR space for tetrahedral and prismatic elements at $P=3$ and $P=6$. Within each subfigure, the HO element is shown on the left and its LOR subdivision on the right. Red points denote the vertices of the piecewise-linear LOR elements. The triangular nodal distributions use electrostatic interior points $\mathbf{\xi_{ele}}$ with boundary points $\mathbf{\xi_{GLL}}$ constrained to the neighbouring edges.}
\label{fig:LOR_simplices}
\end{figure}

The LOR spaces for tetrahedral and prismatic elements require additional connectivity and orientation considerations. Figure~\ref{fig:LOR_simplices} visualises these constructions. Figures~\ref{subfig:Tet_P3} and~\ref{subfig:Tet_P6} show tetrahedral LOR meshes at $P=3$ and $P=6$, respectively.
Figures~\ref{subfig:Prism_P3} and~\ref{subfig:Prism_P6} show prismatic LOR meshes at $P=3$ and $P=6$. The prismatic expansion is constructed as a tensor product of the triangular nodal expansion and a one-dimensional Lagrange expansion.

\subsection{Inner-iteration strategy} 
The strategy for solving the iterative solver residual, interpolated onto the LOR discretisation (Step 6 in Algorithm \ref{alg:DoLORPreconditioner}), or the \textit{inner iteration strategy}, can affect the performance of the LOR preconditioner. A direct solver is the quickest option for small cases and the best approach for running on a single MPI rank. However, for larger problem sizes, the PETSc interface \cite{petsc-web-page} to \nekpp provides a variety of preconditioners and linear solvers, including AMG solvers such as \textit{BoomerAMG} \cite{yang2002boomeramg} and \textit{GAMG}. Native linear solvers and preconditioners in \nekpp can also be applied to the LOR space.

%% file: algos/LOR_process_global.tex
\begin{figure}[H]
    \centering
    \begin{tikzpicture}
        \node[anchor=south west, inner sep=0] (image) at (0,0) {%
            \includegraphics[width=\textwidth]
            {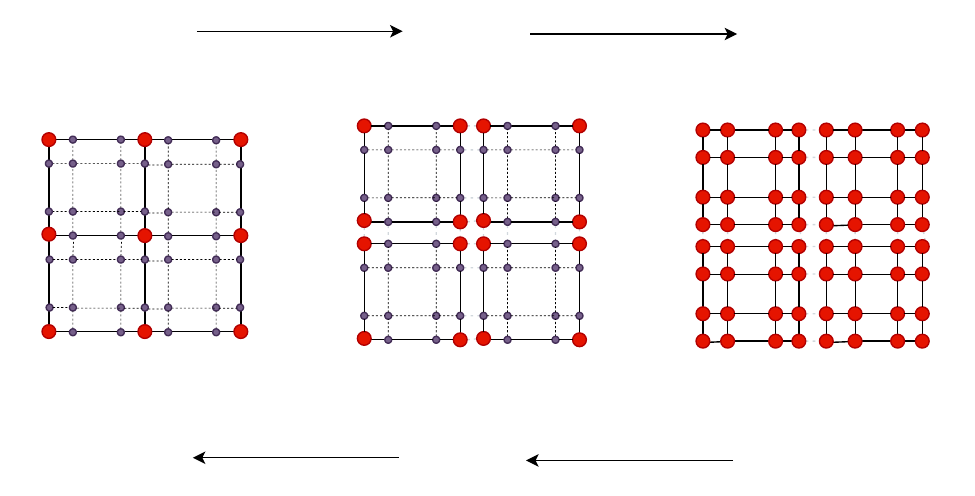}%
        };

        \node[below] at (2.0,5.5)
            {\small Input $\hat{\mathbf{r}}_{H,g}$};

        \node[below] at (10.8,5.6)
            {\small Solve LOR System};

        \node[below] at (10.8,1.9)
            {\small
            $\hat{\mathbf{z}}_{L,g}
            =
            \mathbf{M}_{\mathrm{LOR}}^{-1}
            \hat{\mathbf{r}}_{L,g}$};

        \node[below] at (4.0,6.8)
            {\small Global to Local};

        \node[below] at (4.0,6.2)
            {\small
            $\hat{\mathbf{r}}_{H,g}
            \xrightarrow{\mathcal{A}_H}
            \hat{\mathbf{r}}_{H,l}$};

        \node[below] at (8.4,6.8)
            {\small Change of Basis: HO to LOR};

        \node[below] at (8.4,6.2)
            {\small
            $\hat{\mathbf{r}}_{L,g}
            =
            \mathcal{V}^{-T}\hat{\mathbf{r}}_{H,l}$};

        \node[below] at (4.0,1.1)
            {\small Local to Global};

        \node[below] at (4.0,0.5)
            {\small
            $\hat{\mathbf{z}}_{H,g}
            \xleftarrow{\mathcal{A}_H^{T}}
            \hat{\mathbf{z}}_{H,l}$};

        \node[below] at (8.4,1.1)
            {\small Change of Basis: LOR to HO};

        \node[below] at (8.4,0.5)
            {\small
            $\hat{\mathbf{z}}_{H,l}
            =
            \mathcal{V}^{T}\hat{\mathbf{z}}_{L,g}$};
    \end{tikzpicture}

    \caption{%
        Application of the LOR preconditioner at one preconditioned
        iteration. The global HO residual
        $\hat{\mathbf{r}}_{H,g}$ is mapped to the element-local
        representation by $\mathcal{A}_H$ and transformed to the globally
        assembled LOR representation by $\mathcal{V}^{-T}$. The LOR system
        is then solved, after which the correction is transformed back by
        $\mathcal{V}^{T}$ and assembled in the global HO space by
        $\mathcal{A}_H^{T}$. Solid black lines delineate the local elements,
        while the faint connections between neighbouring elements indicate
        boundary points identified by the global assembly before and after
        LOR refinement. Here, $H$ and $L$ denote the HO and LOR spaces,
        respectively; $g$ and $l$ denote global and element-local
        representations; $\mathcal{V}$ is the generalised Vandermonde; and $\mathbf{M}_{\mathrm{LOR}}^{-1}$
        denotes the LOR solve.
    }
    \label{fig:LOR_process}
\end{figure}

%% file: algos/doLOR.tex
\begin{algorithm}
\caption{Do-Stage: LOR Preconditioning Procedure (per iteration)}
\label{alg:DoLORPreconditioner}
\begin{algorithmic}[1]

\REQUIRE High-order (HO) residual $\hat{\mathbf{r}}_{H,g}$, HO assembly map $\mathcal{A}_H$
\ENSURE Preconditioned residual $\hat{\mathbf{z}}_{H,g}$ in HO coefficient space

\STATE Apply global elementwise multiplicity: $\hat{\mathbf{r}}_{H,g} \gets \mathbf{D}_H^{-1}\hat{\mathbf{r}}_{H,g}$

\STATE Scatter to HO-local coefficients: $\hat{\mathbf{r}}_{H,l} \gets \mathcal{A}_H \hat{\mathbf{r}}_{H,g}$

\STATE Transform to LOR space using the generalized Vandermonde operator: 
$\hat{\mathbf{r}}_{L,g} \gets \mathcal{V}^{-T} \hat{\mathbf{r}}_{H,l}$

\STATE Apply elementwise LOR multiplicity correction: 
$\hat{\mathbf{r}}_{L,g} \gets (\mathbf{D}_L^e)^{-1}\hat{\mathbf{r}}_{L,g}$


\STATE Solve the LOR preconditioner system: 
$\hat{\mathbf{z}}_{L,g} \gets \mathbf{M}_{LOR}^{-1}\hat{\mathbf{r}}_{L,g}$



\STATE Project to HO-local space: 
$\hat{\mathbf{z}}_{H,l} \gets \mathcal{V}^{T}\hat{\mathbf{z}}_{L,g}$

\STATE Assemble to HO-global coefficients (sum of contributions): 
$\hat{\mathbf{z}}_{H,g} \gets \mathcal{A}_H^{T}\hat{\mathbf{z}}_{H,l}$

\RETURN $\hat{\mathbf{z}}_{H,g}$ as the preconditioned residual

\end{algorithmic}
\end{algorithm}

%% file: algos/buildLOR.tex
\begin{algorithm}
\caption{Build stage: LOR preconditioning procedure (one-time operation)}
\label{alg:BuildLORPreconditioner}
\begin{algorithmic}[1]

\REQUIRE High-order mesh graph $\mathcal{M}_H$, HO assembly map $\mathcal{A}_H$, subdivision factor $n_{\text{split}}$
\ENSURE Preconditioner $\mathbf{M}_{\text{LOR}}$, LOR mesh graph $\mathcal{M}_L$, LOR element-wise assembly map $\mathcal{A}_L^e$

\STATE Generate LOR mesh graph: $\mathcal{M}_L \gets \text{Subdivide}(\mathcal{M}_H, n_{\text{split}})$

\STATE Construct LOR assembly and coefficient maps: $\mathcal{A}_L^e, \mathcal{P}_{H\rightarrow L} \gets \text{BuildMaps}(\mathcal{A}_H, \mathcal{M}_L)$

\STATE Assemble LOR system matrix: $\mathbf{H}_L \gets \text{Assemble}(\mathcal{M}_L, \mathcal{A}_L^e)$

\STATE Apply multiplicity scaling: $\mathbf{H}_L \gets \text{Weight}(\mathbf{H}_L, \mathcal{A}_H, \mathcal{A}_L^e)$

\IF{HO basis is Modal}
    \STATE Define interpolation operator: $\mathcal{V} \gets \text{Vandermonde}(\mathcal{M}_H, \mathcal{M}_L)$
    \STATE Construct preconditioner: $\mathbf{M}_{\text{LOR}} \gets \mathcal{V}^{T}\mathbf{H}_L\mathcal{V}$
\ELSE
    \STATE Construct preconditioner (Nodal): $\mathbf{M}_{\text{LOR}} \gets \mathbf{H}_L$
\ENDIF

\end{algorithmic}
\end{algorithm}

%% file: main_conditioning.tex
The primary performance bottleneck in the \IncNS solver is the pressure Poisson system, which serves as the test problem for this study:
\begin{equation}
-\nabla^2 u(\mathbf{x}) = f(\mathbf{x}), \quad \text{with} \quad f(\mathbf{x}) = -d \prod_{i=1}^{d} \sin(x_i).
\label{eq:poisson_combined}
\end{equation}
Unless otherwise specified, we enforce either zero Neumann, zero Dirichlet, or $u_{\partial\Omega} = \prod_{i=1}^{d} \sin(x_i)$ Dirichlet boundary conditions, where $d$ is the spatial dimension. To characterize the spectral properties of the preconditioned system, we compute the iterative condition number:
\begin{equation}
\kappa(\mathbf{M}^{-1}\mathbf{H}) = \frac{|\lambda_{\max}|}{|\lambda_{\min}|}.
\end{equation}
The extreme eigenvalues $\lambda$ are estimated via a standard QR-based eigensolver \cite{lapack99} to assess the effectiveness of each preconditioning strategy $\mathbf{M}$.

\subsection{Conditioning properties}
\label{subsec:LOR_validation}

\subsubsection{2D test cases}
A single-element 2D unit-square case is the first test case of this study. The Poisson problem in Eqn.~\eqref{eq:poisson_combined} is solved for $d = 2$. Table ~\ref{tab:quad1_conditionnumber} has the condition number results for a p-convergence study of the LOR preconditioner with zero Dirichlet boundary conditions (BCs) on all domain boundaries for the unit-square.

The results of a $\mathds{Q}1$ discretisation using a nodal expansion basis are in close agreement with those of a similar previous implementation by \cite{PedroFischer2019}. It follows the condition number convergence with increasing problem size ($\kappa \approx \pi^2/4$)  for Dirichlet cases mentioned by \cite{Canuto2009}. The $\mathds{P}1$ discretisation with a modal expansion basis also converges to the same value. However, that is not the case for the $\mathds{Q}1$ discretisation, which converges to a higher value. Some differences in modal expansion basis results are expected as the change of basis operation is done via a generalised Vandermonde matrix projection, which may introduce some interpolation errors.

\begin{table}[H]
\centering
\caption{Comparison of the condition number ($\kappa$) for the 2D unit-square case, single spectral element with Dirichlet BCs on all boundaries}
\label{tab:quad1_conditionnumber}
\begin{tabular}{ccccc}
\hline
\multicolumn{1}{l}{\textbf{\textit{P}}} &
  \multicolumn{1}{l}{\textbf{\cite{PedroFischer2019}}} &
  \multicolumn{1}{l}{\textbf{Nodal, $\mathds{Q}1$}} &
  \multicolumn{1}{l}{\textbf{Modal, $\mathds{Q}1$}} &
  \multicolumn{1}{l}{\textbf{Modal , $\mathds{P}1$}} \\ \hline
\textbf{2}  & 1     & 1     & 1     & 1     \\
\textbf{4}  & 1.554 & 1.555 & 2.466 & 1.392 \\
\textbf{6}  & 1.803 & 1.804 & 3.252 & 1.692 \\
\textbf{8}  & 1.954 & 1.945 & 3.782 & 1.869 \\
\textbf{10} & 2.037 & 2.037 & 4.418 & 1.982 \\
\textbf{12} & 2.101 & 2.101 & 4.877 & 2.087 \\
\textbf{14} & 2.148 & 2.149 & 5.219 & 2.192 
\end{tabular}
\end{table}

\begin{figure}[H]
\centering
\begin{subfigure}{.2\textwidth}
  \centering
  \raisebox{8mm}{\includegraphics[width=\linewidth]{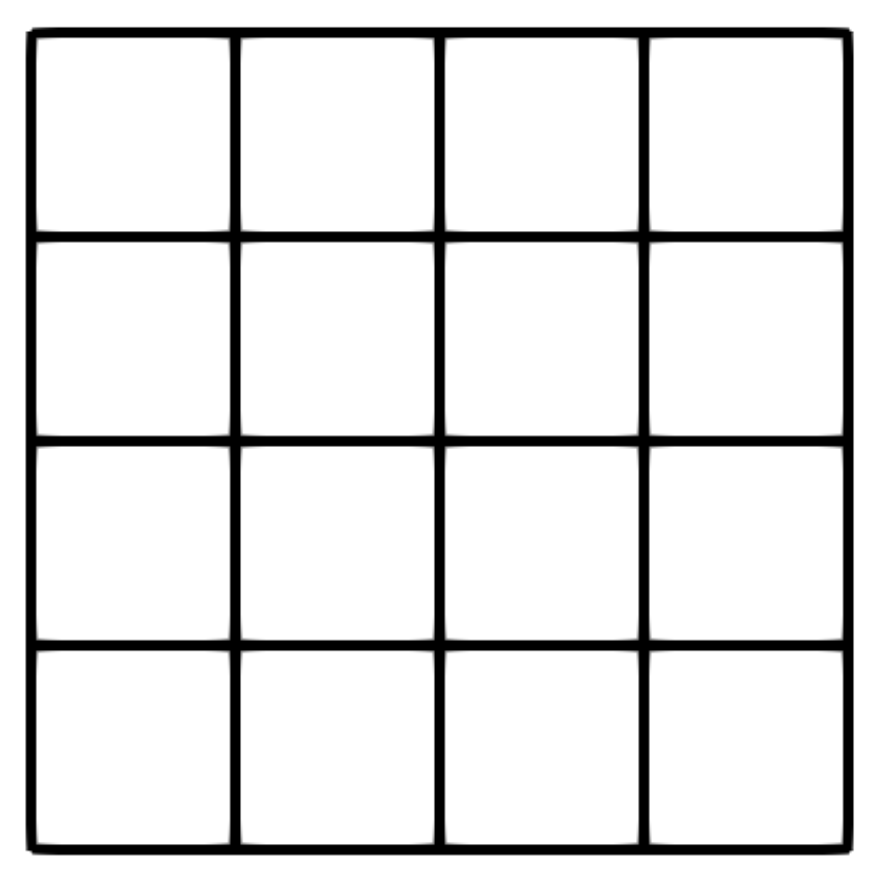}}
  \caption{Quad16 mesh}
  \label{subfig:Quad16}
\end{subfigure}%
\hfill
\begin{subfigure}{.2\textwidth}
  \centering
  \raisebox{8mm}{\includegraphics[width=\linewidth]{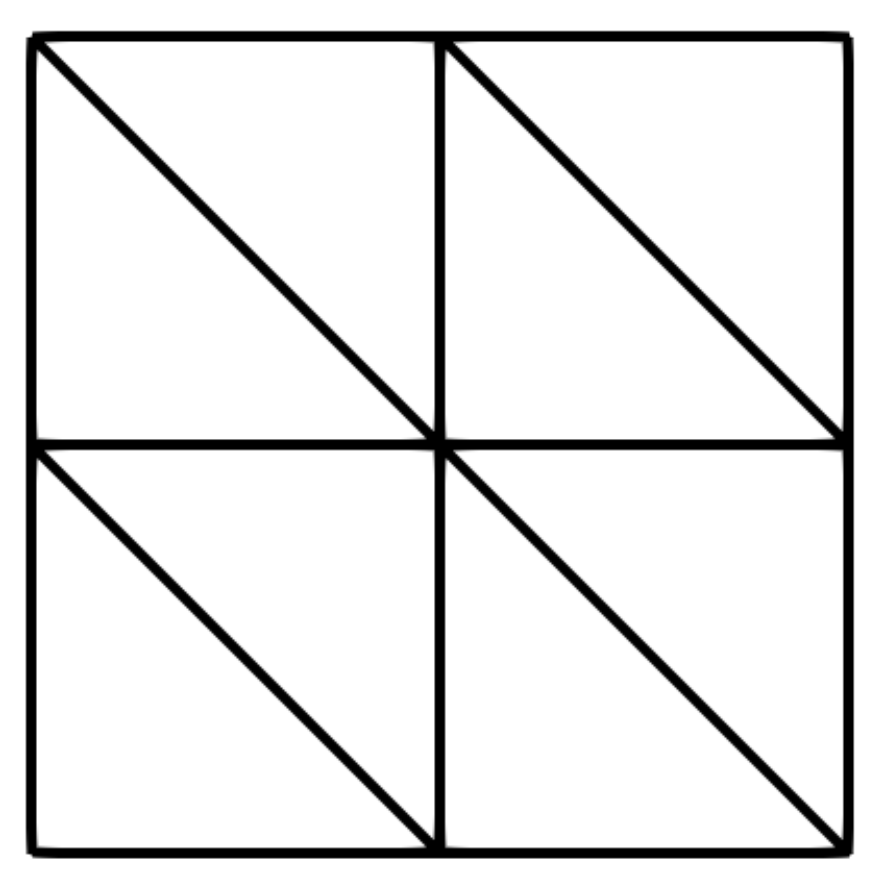}}
  \caption{Tri8 mesh}
  \label{subfig:Tri8}
\end{subfigure}%
\hfill
\begin{subfigure}{.5\textwidth}
  \centering
  \includegraphics[width=\linewidth]{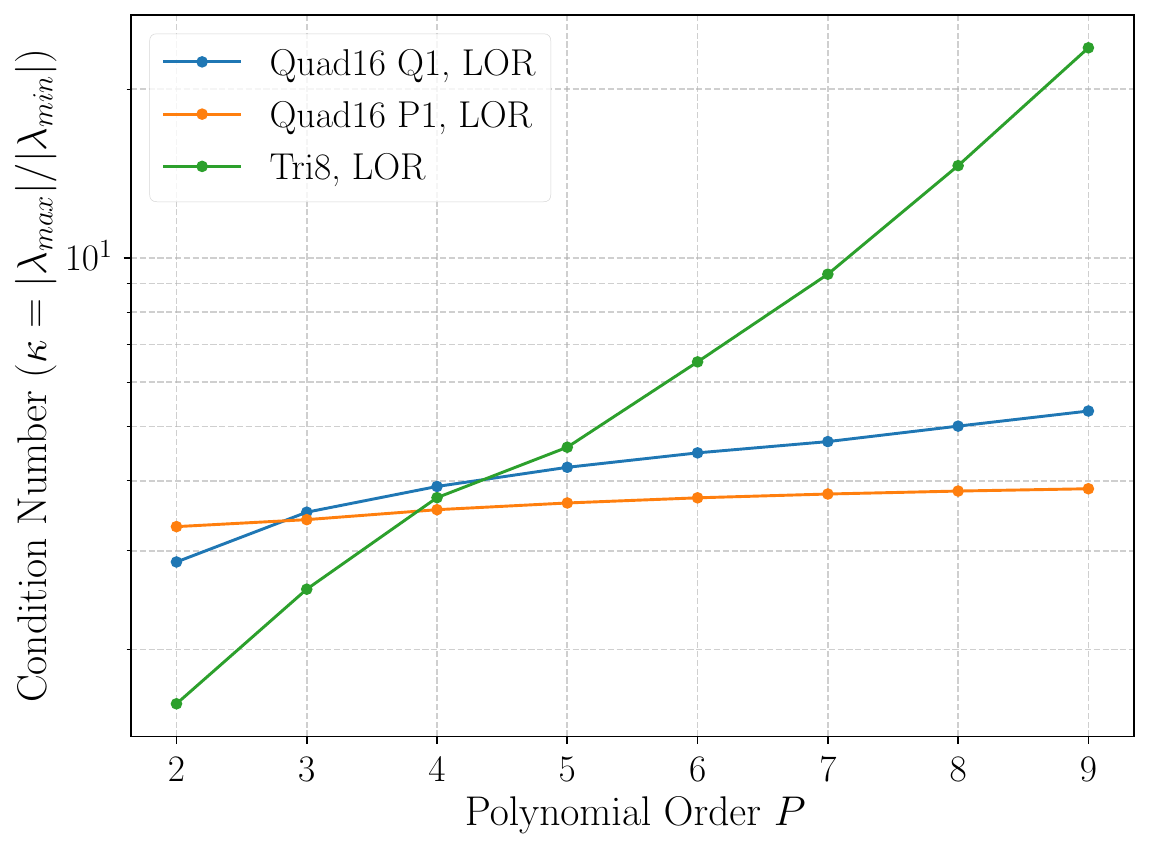}
  \caption{$\kappa$ with \textit{P} for the two meshes}
  \label{subfig:condnum_P_quad}
\end{subfigure}
\caption{Comparison of the condition number ($\kappa$) on applying the LOR preconditioner for the 2D unit-square case with zero Dirichlet BCs on all boundaries. Two different meshes are investigated to see the effect of different types of higher-order elements.}
\label{fig:quad}
\end{figure}


A condition number analysis for increasing polynomial orders is conducted on 2D meshes (Fig.~\ref{fig:quad}) to study the influence of element type on conditioning in 2D geometries. The Poisson problem in Eqn.~\eqref{eq:poisson_combined} is solved for $d = 2$, with zero Dirichlet boundary conditions applied on all edges. From this point onward, only modal expansions are considered, as they represent the most common use case within the \nekpp framework. The triangular elements (Tri8 mesh, Fig.~\ref{subfig:Tri8}) show an exponentially increasing condition number trend, and exceed the quadrilateral elements (Quad16 mesh, Fig.~\ref{subfig:Quad16}) at $P\geq 5$. For the Quad16 case, using the $\mathds{P}1$ element type yields marginally improved conditioning. However, for the polynomial range of $P \leq 6$, the condition number is controlled.

\begin{figure}[H]
    \centering
    \includegraphics[width=0.7\textwidth]{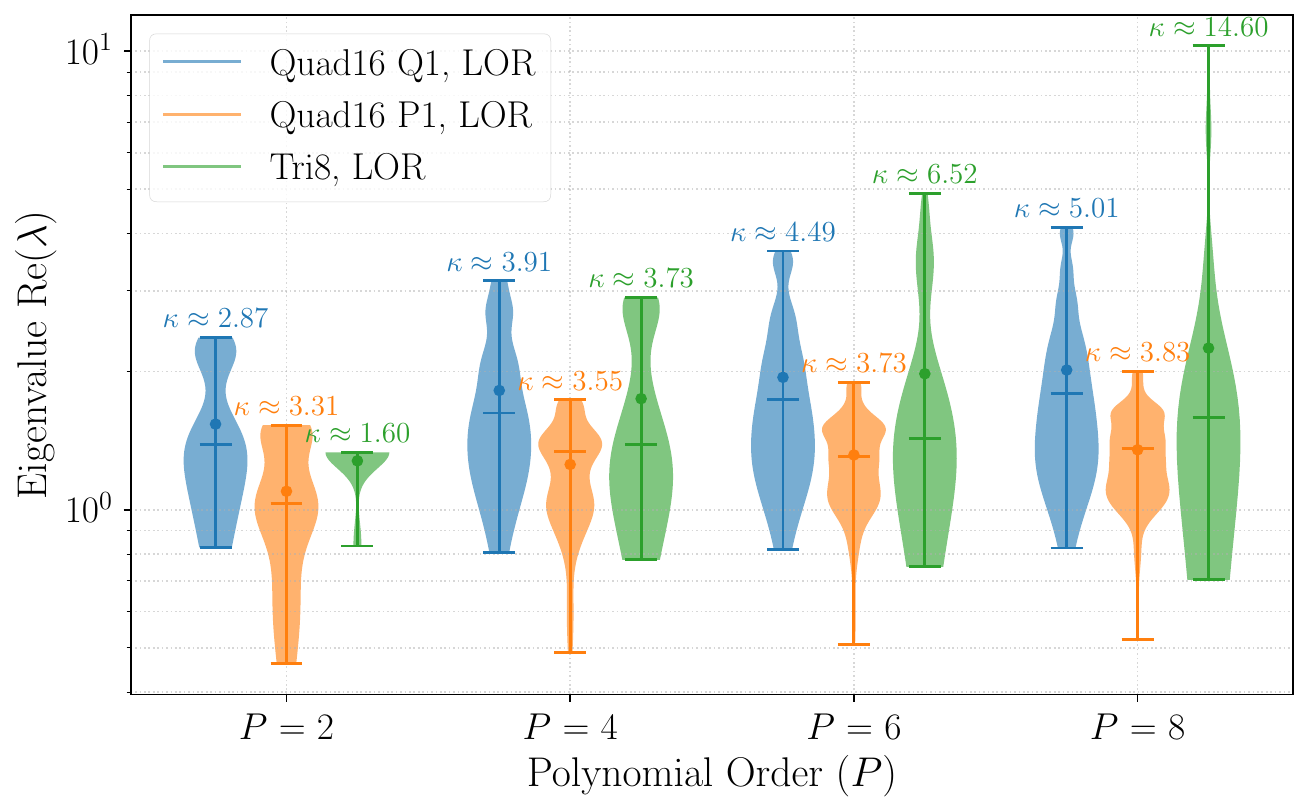} 
    \caption{Eigenvalue spread for the 2D meshes in Fig.~\ref{fig:quad} on applying the LOR preconditioner for different polynomial orders. $\kappa$ = iterative condition number of the preconditioned matrix, or $\kappa(\mathbf{M}^{-1}\mathbf{H})$. All eigenvalues are real numbers, from the Poisson equation's symmetric positive definite (SPD) matrix. The dot on the violin plot indicates the mean value, the horizontal lines are the median and the extrema, and the thickness is the density of the distribution.}
\label{fig:EV_spread_2D}
\end{figure}

The corresponding eigenvalue distributions of the LOR-preconditioned operator (Fig.~\ref{fig:EV_spread_2D}) further clarify this behaviour. The eigenvalues are expected to be real numbers from the Poisson equation's symmetric positive definite (SPD) nature. The minimum eigenvalue stays constant across all cases for increasing polynomial orders. The maximum eigenvalue is controlled for the mesh with quadrilateral elements but increases for the triangular element mesh. The low frequency component is adequately addressed for both element types, but the high frequency component may become difficult to condition as the increased anisotropy of the LOR space at high $P$ is more severe for triangular discretisation.

\subsubsection{3D test cases}
\begin{figure}[H]
\centering
\begin{subfigure}{.25\textwidth}
  \centering
  \includegraphics[width=\linewidth]{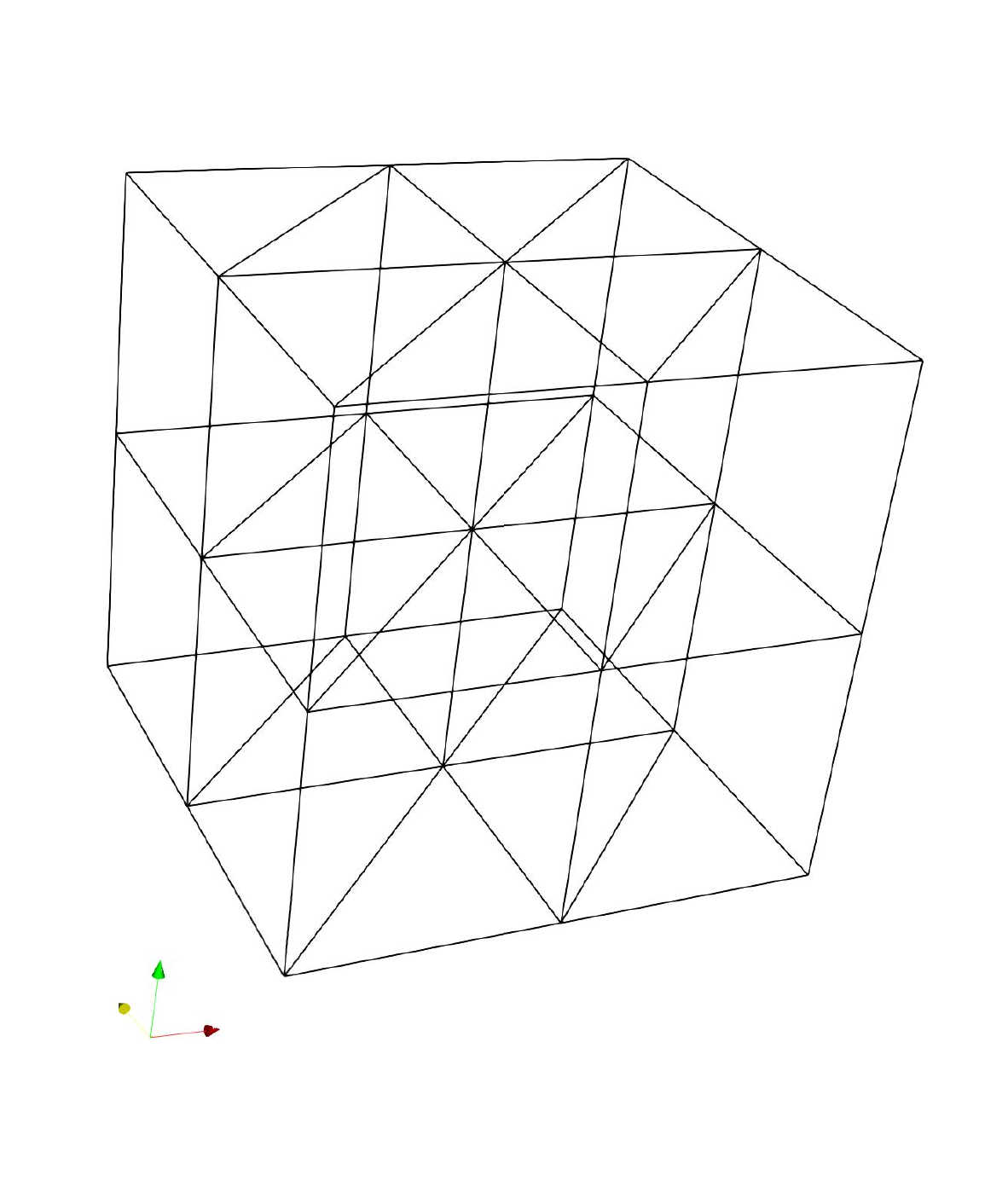}
  \caption{Cube-prism mesh}
  \label{subfig:cube-prism}
\end{subfigure}%
\begin{subfigure}{.25\textwidth}
  \centering
  \includegraphics[width=\linewidth]{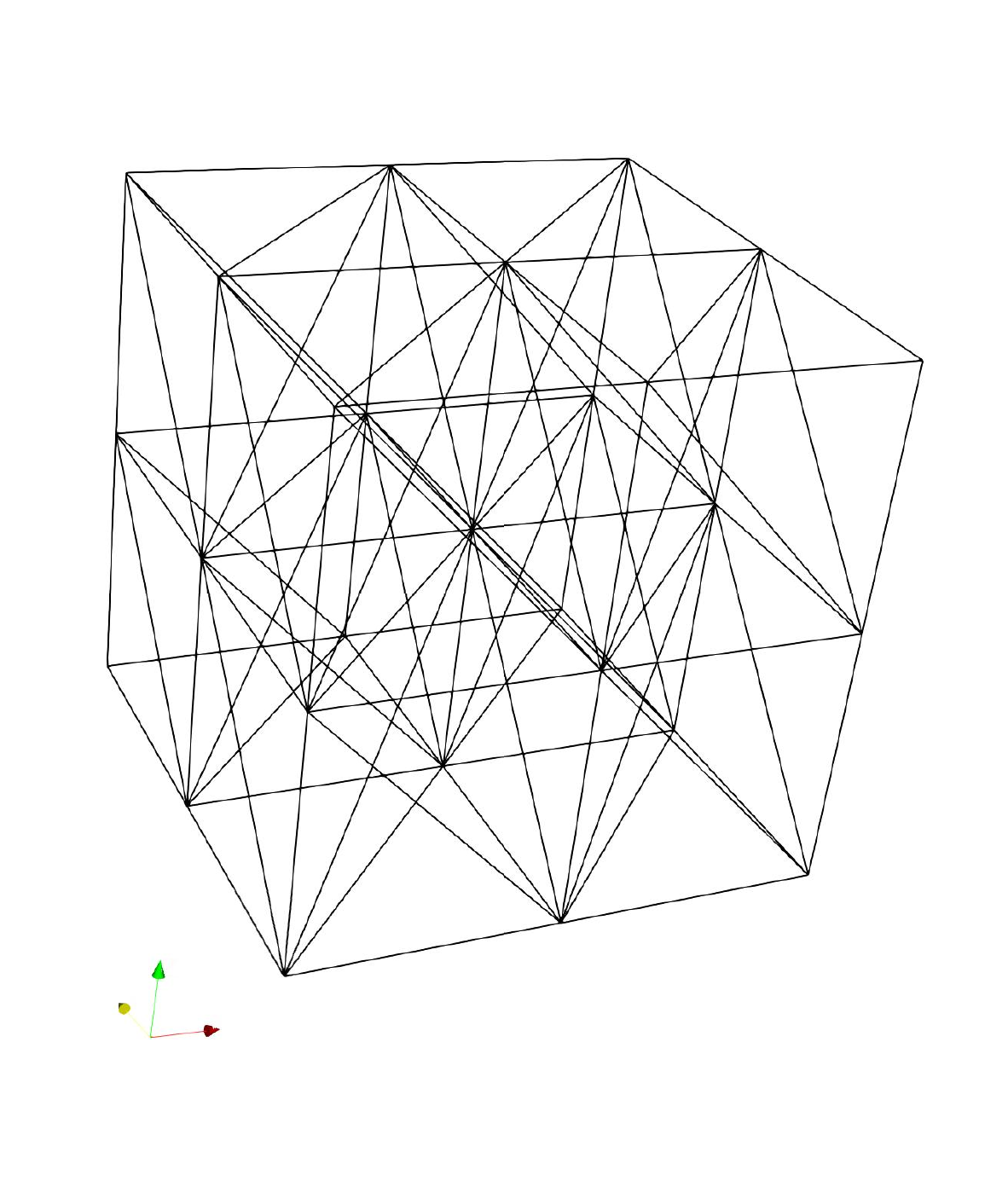}
  \caption{Cube-tet mesh}
  \label{subfig:cube-tet}
\end{subfigure}
\begin{subfigure}{.45\textwidth}
  \centering
  \includegraphics[width=\linewidth]{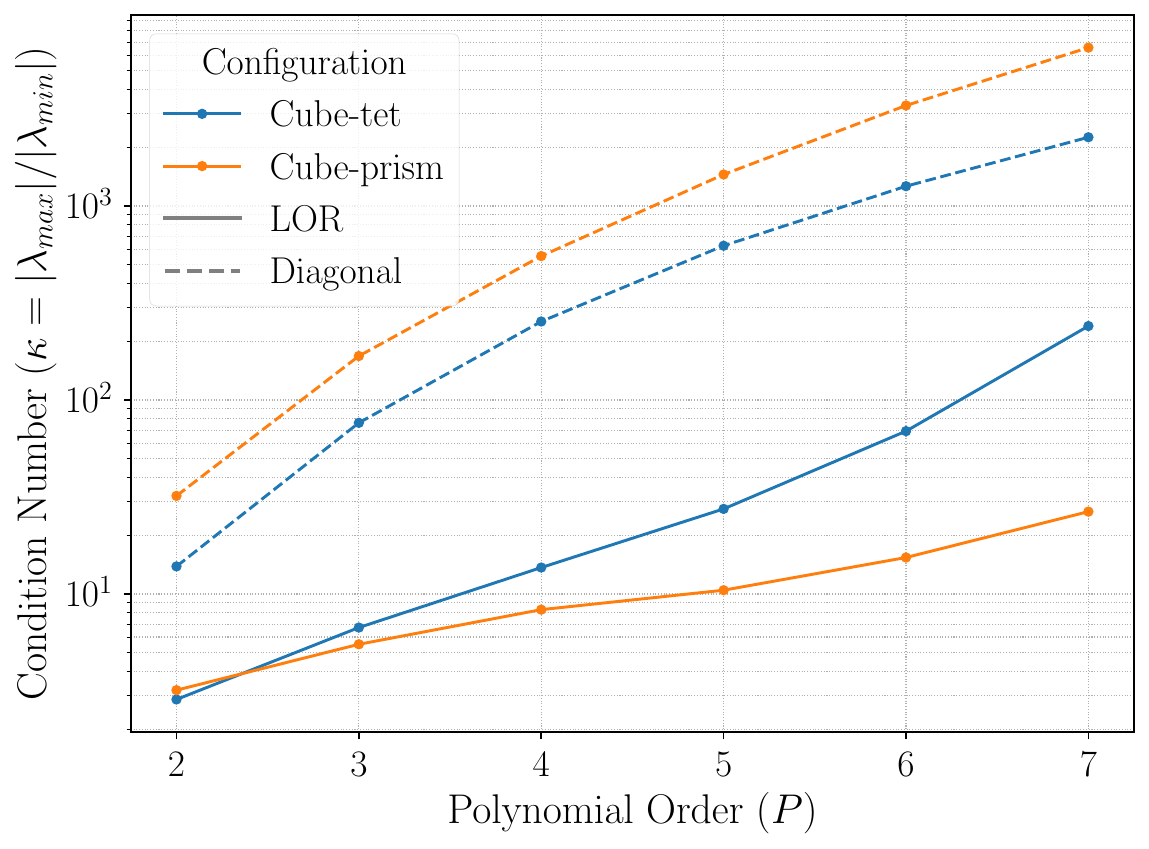}
  \caption{$\kappa$ with \textit{P} for different preconditioners}
  \label{subfig:condnum_P}
\end{subfigure}
\caption{Comparison of the condition number ($\kappa$) on applying the LOR preconditioner for the cube case with Dirichlet BCs on all boundaries. Two different meshes are investigated to see the effect of different types of higher-order elements.}
\label{fig:condnum_cube}
\end{figure}

A unit cube in 3D, as seen in Fig.~\ref{fig:condnum_cube}, discretised using tetrahedra and prisms, is the first 3D test case of this study. The LOR preconditioner is evaluated on meshes containing only prisms (Fig.~\ref{subfig:cube-prism}) and only tetrahedra (Fig.~\ref{subfig:cube-tet}) to understand the conditioning properties for these element types. The Poisson problem in Eqn.~\eqref{eq:poisson_combined} is solved for $d = 3$ and Dirichlet BCs $u_{\partial\Omega} = \prod_{i=1}^{d} \sin(x_i)$ on all the domain boundaries. A modal HO basis function is chosen.

Fig.~\ref{subfig:condnum_P} compares the condition number for a prism mesh (Cube-prism) to the tetrahedron mesh (Cube-tet) for increasing polynomial order. The LOR preconditioner consistently outperforms the Diagonal preconditioner, with condition numbers lower by more than an order of magnitude as the polynomial degree increases. While the condition number still grows with $P$, the trend remains controlled up to $P \leq 6$. Notably, in the range $P \in [3,6]$, the significant gap in conditioning between LOR and Diagonal highlights the effectiveness of the LOR approach, making it a promising candidate for use in complex 3D geometries.

\begin{figure}[H]
    \centering
    \includegraphics[width=0.6\textwidth]{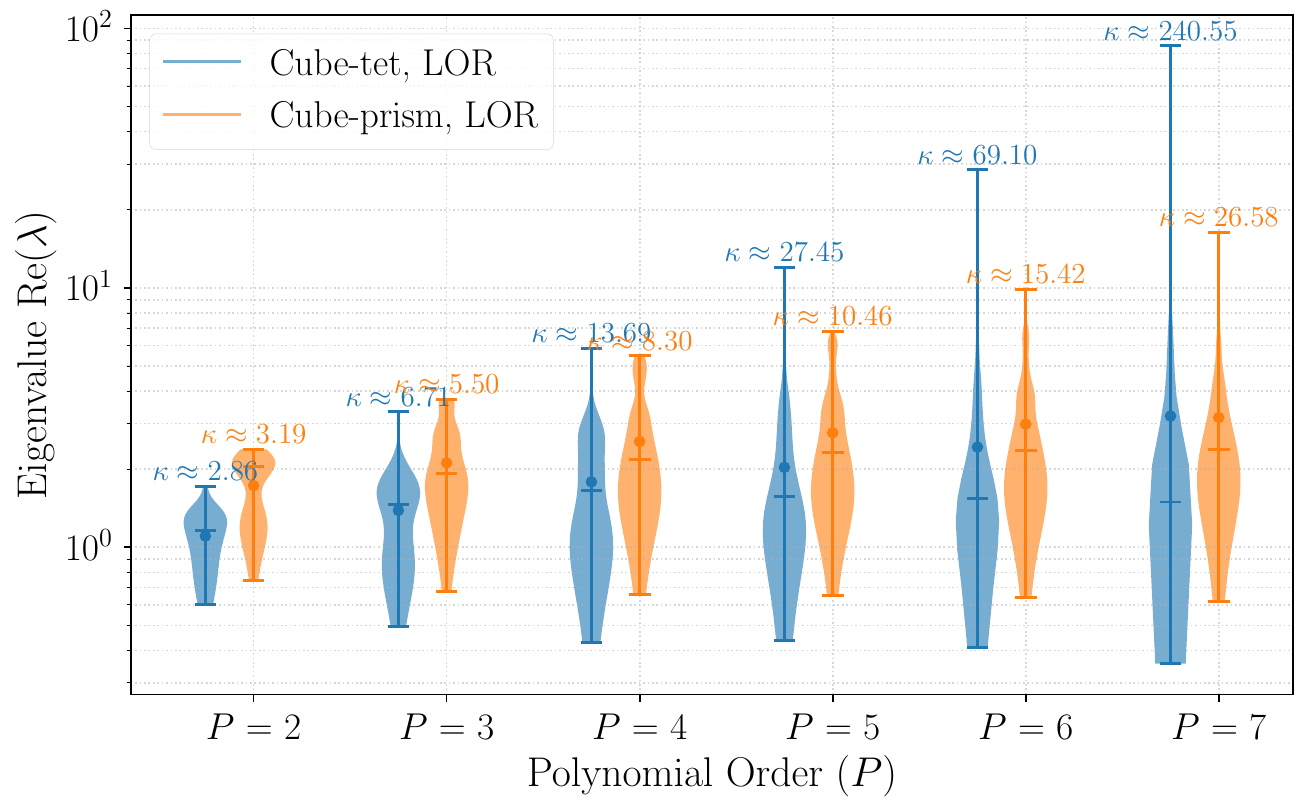} 
    \caption{Eigenvalue spread for the 3D meshes in Fig.~\ref{fig:condnum_cube} on applying the LOR preconditioner for different polynomial orders. $\kappa$ = iterative condition number of the preconditioned matrix, or $\kappa(\mathbf{M}^{-1}\mathbf{H})$. All eigenvalues are real numbers, from the Poisson equation's symmetric positive definite (SPD) matrix. The dot on the violin plot indicates the mean value, the horizontal lines denote the median and the extrema, and the thickness of the violin indicates the distribution density.}
    \label{fig:EV_spread_3D}
\end{figure}

Fig.~\ref{fig:EV_spread_3D} presents the eigenvalue distribution for the two meshes across various polynomial orders. For the prism elements, the minimum eigenvalue remains approximately constant while the maximum eigenvalue increases as the polynomial order increases. For tetrahedral elements, the minimum eigenvalue decreases, and the maximum eigenvalue increases sharply, indicating worse conditioning. For both cases, the density of the clustering remains the same, and only the maximum eigenvalue shoots up with the polynomial order increase.

%% file: main_iterative.tex
\subsection{Iterative performance}
\label{sec:LOR_iterative_performance}
To further assess the iterative properties of the LOR preconditioner, the performance of the preconditioned GMRES solution is evaluated across different geometries. The test cases evaluated include various mixed-element meshes in 2D and 3D. Comparisons are made with various iterative solver configurations in \nekpp. All the tests in this section, which use an iterative solver, have been run until a convergence tolerance $10^{-4}$ is reached. The solver used most frequently is GMRES, with a restart after the GMRES projection reaches 100 vectors. The \textit{BoomerAMG} settings for the AMG solver are specified in Appendix \ref{App:petscrc_settings}.

\subsubsection{2D NACA aerofoil}
\begin{figure}[H]
\centering
\begin{subfigure}{0.45\textwidth}
  \centering
  \includegraphics[width=\linewidth]{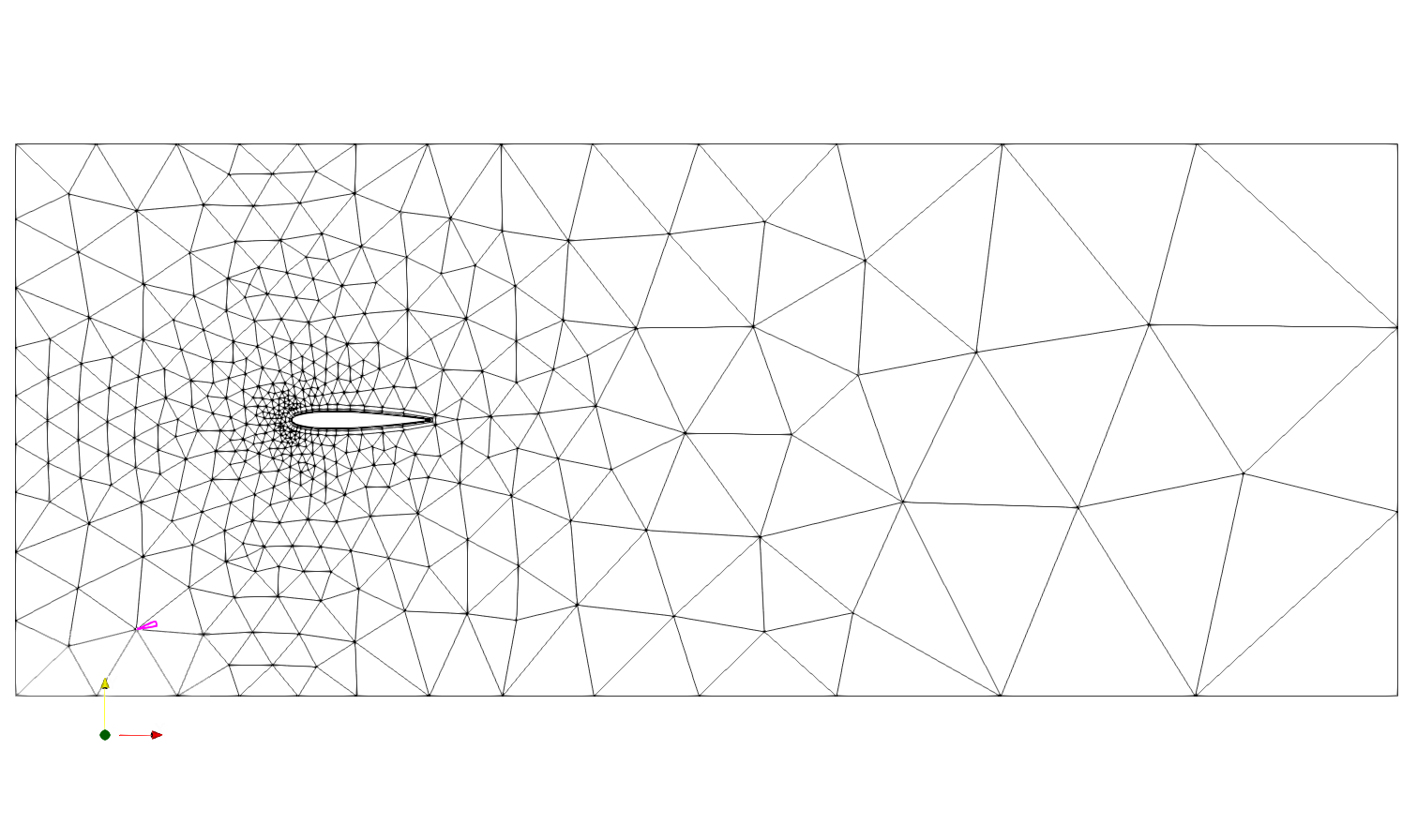}
  \caption{Full Domain, HO mesh ($N_{ele} = 876$)}
  \label{subfig:naca_FullHO}
\end{subfigure}%
\hfill
\begin{subfigure}{.45\textwidth}
  \centering
  \includegraphics[width=\linewidth]{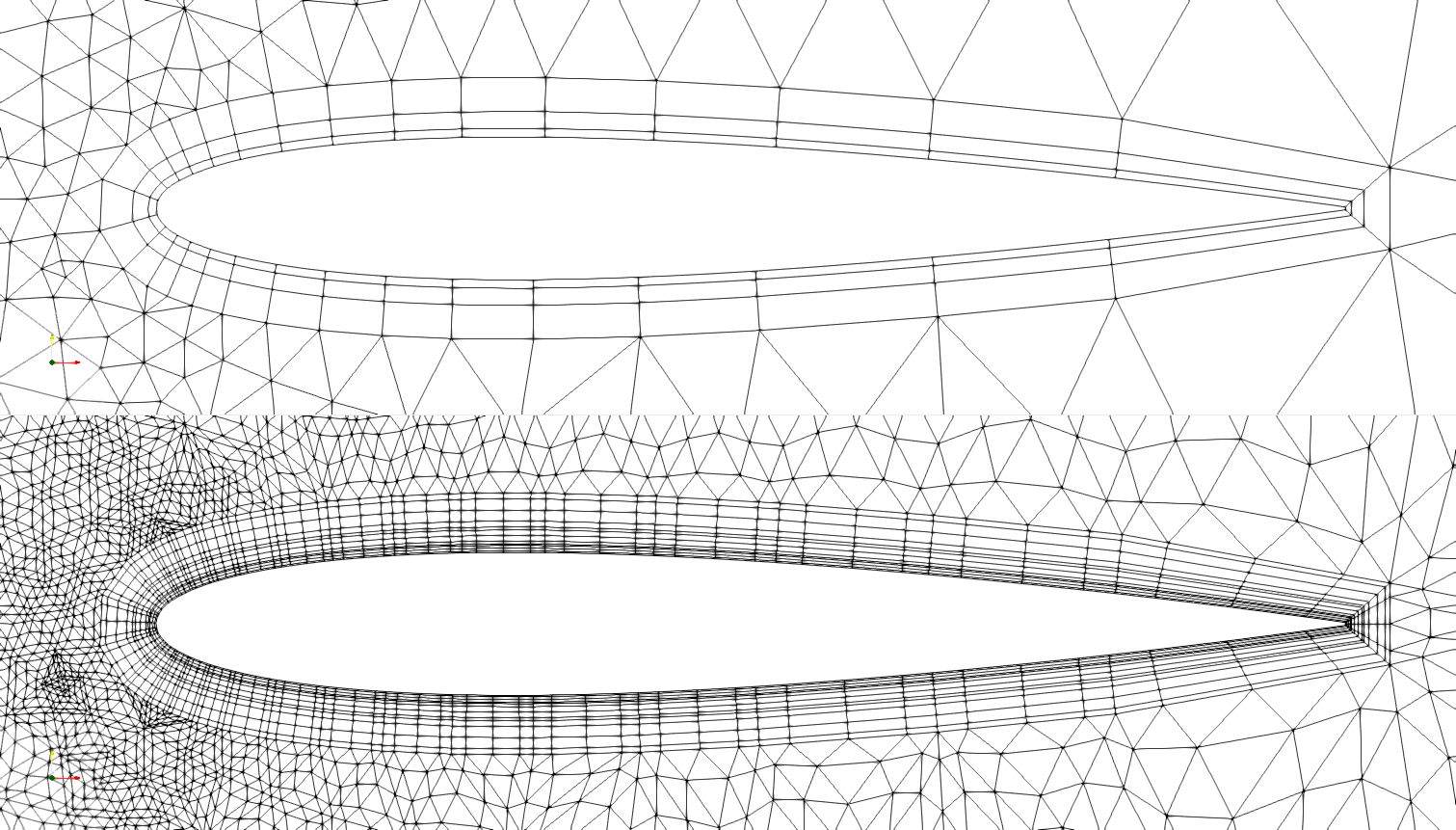}
  \caption{HO (above) vs LOR (below) at P=4}
  \label{subfig:naca_LORvsHO}
\end{subfigure}
\caption{Domain for the 2D NACA0012 airfoil test case. The mesh is a combination of quadrilaterals and triangles, with three boundary layers along the airfoil profile consisting of quadrilateral elements, and the remainder of the domain is filled with triangles. The LOR mesh at \textit{P}=4 is shown on the right, with a GLL spacing within a given element.}
\label{fig:naca0012_2D}
\end{figure}

The first test case is a flow past the 2D NACA aerofoil case, which is the first mixed-element test case considered in this study (Fig.~\ref{fig:naca0012_2D}). The Poisson problem in Eqn.~\eqref{eq:poisson_combined} is solved for $d = 2$, with Dirichlet BCs $u_{\partial\Omega} = \prod_{i=1}^{d} \sin(x_i)$ on all domain boundaries. The mesh for this test case has a refinement region at the leading edge, and elongated elements in the boundary layers, providing a wide range of element aspect ratios to assess LOR preconditioner performance. This test case compares the performance of the LOR preconditioner with other preconditioning strategies available in \nekpp. The quadrilaterals are chosen to be $\mathds{Q}1$ type in the LOR space. 

Table ~\ref{subtab:Iters_NACA0012} compares the various solution strategies for the test case for solving the Poisson problem at \textit{P} = 4 (7916 DOF). LOR preconditioner performs the best in terms of iterations to the solution. In the polynomial order sweep for this test case in Table ~\ref{subtab:polysweep_NACA0012}, the bounded iteration behaviour with increasing problem size is not seen anymore. The loss of bounded behaviour is attributed to triangle elements in the mixed mesh, which do not have an optimum node set in the LOR space that achieves this behaviour. However, the iteration count is controlled, which is encouraging, especially since the LOR preconditioner is intended for moderate polynomial orders at large problem sizes.

\begin{table}[H]
    \centering
    \caption{Combined iteration count results for the 2D NACA0012 case. The inner loop has been solved till a tolerance of $10^{-7}$ is reached. The inner loop for LOR preconditioner is solved using a direct solver after preconditioning.}
    \label{tab:combined_NACA0012_iters}
    
    \begin{subtable}[t]{0.65\textwidth} 
        \centering
        \caption{Comparison of solver strategies at \textit{P} = 4 (7916 DOF).  }
        \label{subtab:Iters_NACA0012}
        \begin{tabular}{ccc}
            \hline
            \multicolumn{1}{l}{\textbf{Solver}} & \multicolumn{1}{l}{\textbf{Preconditioner}} & \multicolumn{1}{l}{\textbf{Iterations}} \\ \hline
            GMRES & Null & 364 \\
            GMRES & Diagonal & 94 \\
            GMRES & Block & 91 \\
            \textit{BoomerAMG} & - & 72 \\
            PETSc-GMRES & \textit{BoomerAMG} & 21 \\
            CG & LOR & 16 \\
            GMRES & LOR & 14 \\
        \end{tabular}
    \end{subtable}
    \hfill
    \begin{subtable}[t]{0.3\textwidth} 
        \centering
        \caption{Polynomial-order sweep - LOR preconditioner.}
        \label{subtab:polysweep_NACA0012}
        \begin{tabular}{ccc}
            \hline
            \textbf{$P$} & \textbf{GMRES} & \textbf{CG} \\ \hline
            2 & 12 & 14 \\
            3 & 13 & 20 \\
            4 & 14 & 16 \\
            5 & 14 & 16 \\
            6 & 17 & 18 \\
            7 & 18 & 20 \\
        \end{tabular}
    \end{subtable}
\end{table}

\subsubsection{3D Hemisphere}
\begin{figure}[H]
\centering
\begin{subfigure}{.48\textwidth}
  \centering
  \includegraphics[width=\linewidth]{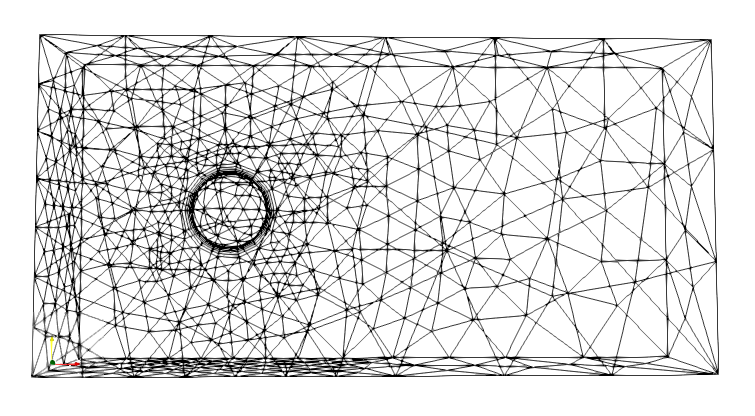}
\end{subfigure}
\begin{subfigure}{.5\textwidth}
  \centering
  \includegraphics[width=\linewidth]{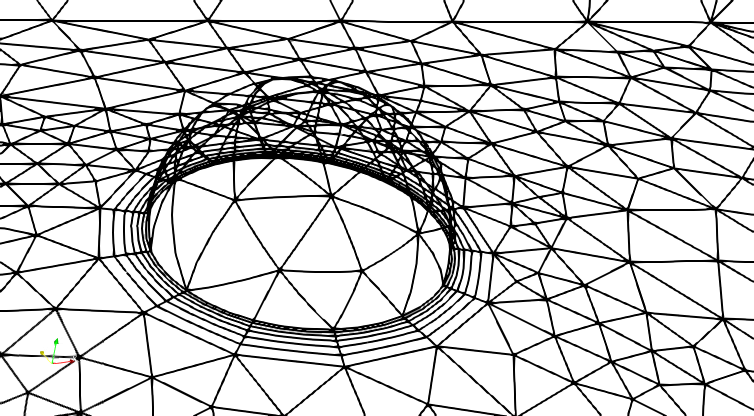}
\end{subfigure}
\caption{Meshes for the 3D Hemisphere case ($N_{ele} = 10,602$). The meshes were generated using the MCF file capability in the mesh generation software \textit{NekMesh} \cite{NekMeshGREEN2024}. Right: the isometric view of the surface mesh inside the domain, Left: the bottom view of the domain. Six prism layers are added to the surface of the hemisphere to create a mesh of Tetrahedra and Prisms.}
\label{fig:sphere}
\end{figure}

To further investigate the iterative properties of the LOR preconditioner in 3D, the next geometry is a flow past a 3D Hemisphere, as seen in Fig.~\ref{fig:sphere}. The three-dimensional ($d=3$) Poisson problem~\eqref{eq:poisson_combined} is evaluated under two configurations: pure Dirichlet conditions ($u_{\partial\Omega} = \prod_{i=1}^{d} \sin(x_i)$) on all boundaries, $\partial\Omega$, and a mixed setup consisting of an outlet Dirichlet boundary with zero Neumann conditions elsewhere. 
The mesh for this geometry is seen in Fig.~\ref{fig:sphere}, where six additional prism layers on the surface of the hemisphere are added. This mesh is a step towards meshes typical of aerodynamics simulations, where boundary layers are added to better resolve the near-wall flow.

\begin{figure}[H]
\centering
  \includegraphics[width=0.9\linewidth]{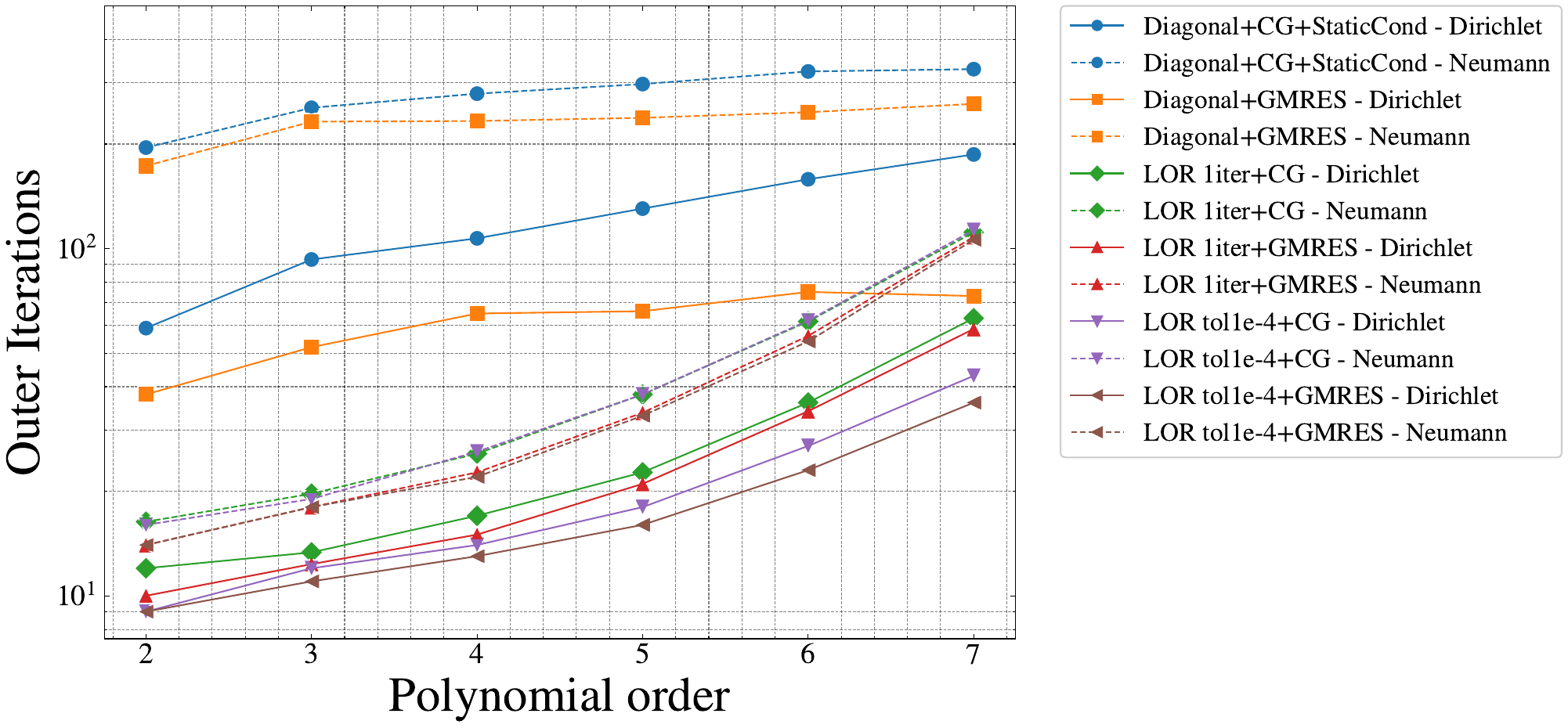}
  \label{subfig:iters_sphereBL}
\caption{Outer iterations to solve the Poisson equation for the 3D Hemisphere case with BL (tetrahedra and prisms) -  polynomial order sweep. The iteration results are averaged over three runs on 10 MPI ranks on a single computational node. StaticCond = Static Condensation. LOR/Diagonal: Preconditioner. CG/GMRES: Iterative solver. Dirichlet: Dirichlet BC on all boundaries. Neumann: Zero Neumann BCs on all domain boundaries apart from the outlet, which is Dirichlet. 1iter = 1 V-cycle of AMG for LOR solution. tol1e-4 = V-cycles for relative tolerance drop of $10^{-4}$}
\label{fig:iters_sphere}
\end{figure}

Fig.~\ref{fig:iters_sphere} shows the results for the outer iterations needed to solve the Poisson equation for the 3D hemisphere case with increasing polynomial order. A comparison is made with the default iterative solver configuration for the pressure system in the \IncNS solver for large 3D cases, also referred to as the ``Baseline", moving forward. The default setup applies a static condensation to the system matrix \cite{karniadakis2005spectral} and uses a conjugate gradient solver with the Diagonal preconditioner (Diagonal+CG+StaticCond in Fig.~\ref{fig:iters_sphere}).

The iterative performance of the LOR preconditioner is consistently best across the two meshes and boundary condition choices for the polynomial range of $P \leq 6$. For higher polynomial orders, the Diagonal preconditioner with GMRES performs better than LOR only for the Dirichlet BCs. Interestingly, LOR consistently outperforms the Baseline configuration.

Two LOR inner iteration strategies are also tested, where at every inner solver iteration, either 1 AMG V-cycle is applied (\textit{1iter}), or V-cycles are applied until a drop of $10^{-4}$ is seen in the residual (\textit{tol1e-4}). The \textit{tol1e-4} approach yields a lower outer iteration count because the inner solution is solved to a lower residual, enabled by multiple V-cycles. The \textit{1iter} approach also yielded variable iteration counts across the three runs as LOR system is not solved to a low enough residual. However, the \textit{1iter} approach is not considerably far off iterations-wise from the \textit{tol1e-4}. A broader comparison of inner-solver cost and convergence is reported in \cite{khurana_2026}.

%% file: main_scaling.tex
\begin{figure}[H]
\centering
\includegraphics[width=\textwidth]{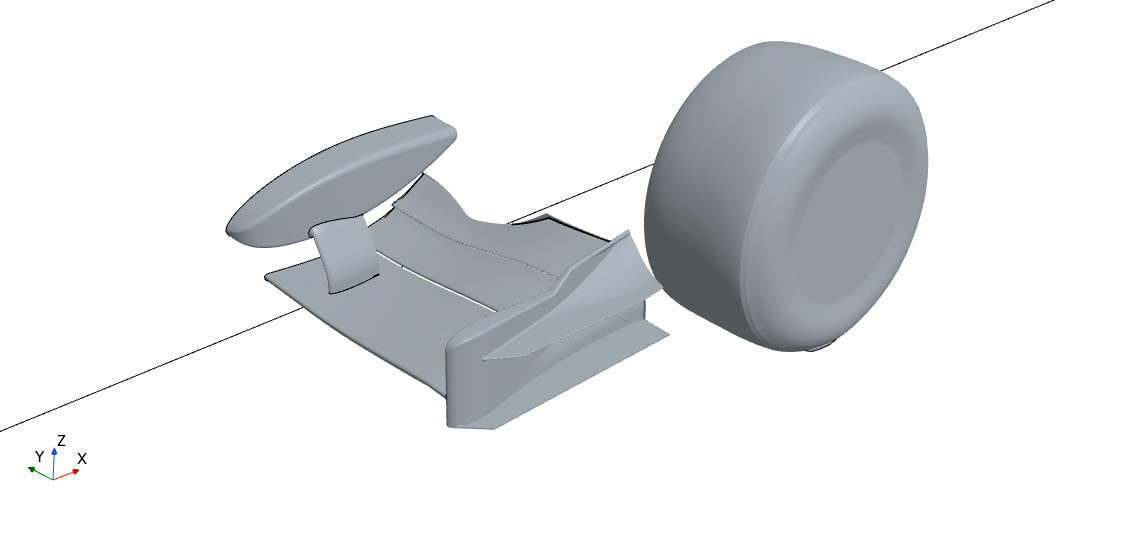}
\caption{Geometry for the IFW with Wheel ($2.87 \times10^6$ elements). The higher-order meshes comprise a mixture of tetrahedra and prisms, with eight boundary layers along the airfoil elements and five boundary layers on the floor profile, which contains prismatic elements. The rest of the domain is filled with tetrahedron elements.}
\label{fig:testcases_IFW}
\end{figure}

The strong-scaling study in this section is carried out on the ARCHER2 cluster, the UK's national supercomputing service. The hardware specification is given in Appendix~\ref{App:hardware}, and the software environment in Appendix~\ref{App:software}. The industrial test case is IFW with a wheel downstream (IFW-W, Fig.~\ref{fig:testcases_IFW}), whose higher-order mesh is a mixture of tetrahedra and prisms.

\subsection{Strong scaling performance}
\label{subsec:wifw_strong_scaling}
\begin{figure}[H]
\centering

\begin{subfigure}[t]{0.48\textwidth}
  \centering
  \includegraphics[width=\linewidth]{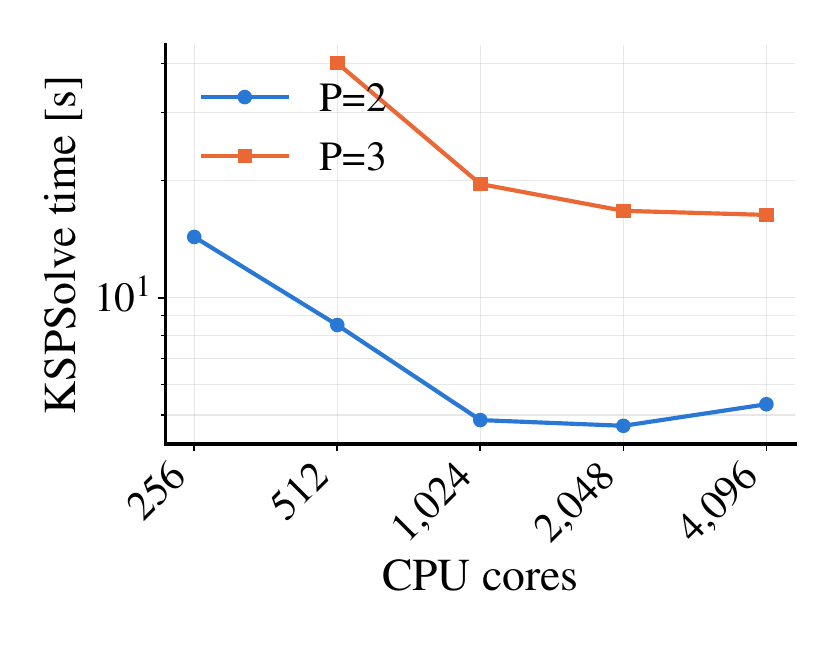}
  \caption{Solve time to solution}
  \label{fig:wifw_scaling_time}
\end{subfigure}
\hfill
\begin{subfigure}[t]{0.48\textwidth}
  \centering
  \includegraphics[width=\linewidth]{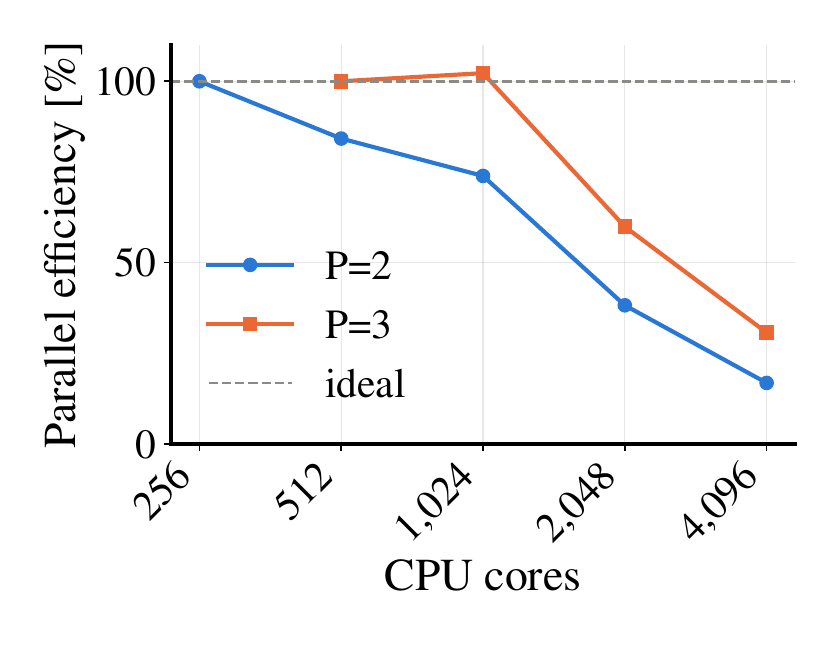}
  \caption{Parallel efficiency}
  \label{fig:wifw_scaling_efficiency}
\end{subfigure}


\begin{subfigure}[t]{0.48\textwidth}
  \centering
  \includegraphics[width=\linewidth]{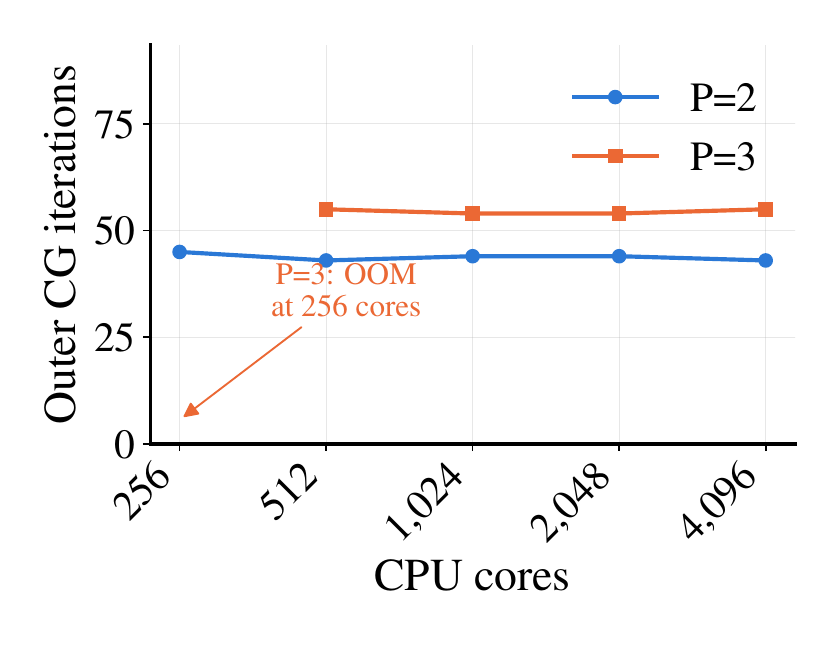}
  \caption{Outer CG iteration count}
  \label{fig:wifw_scaling_iterations}
\end{subfigure}
\hfill
\begin{subfigure}[t]{0.48\textwidth}
  \centering
  \includegraphics[width=\linewidth]{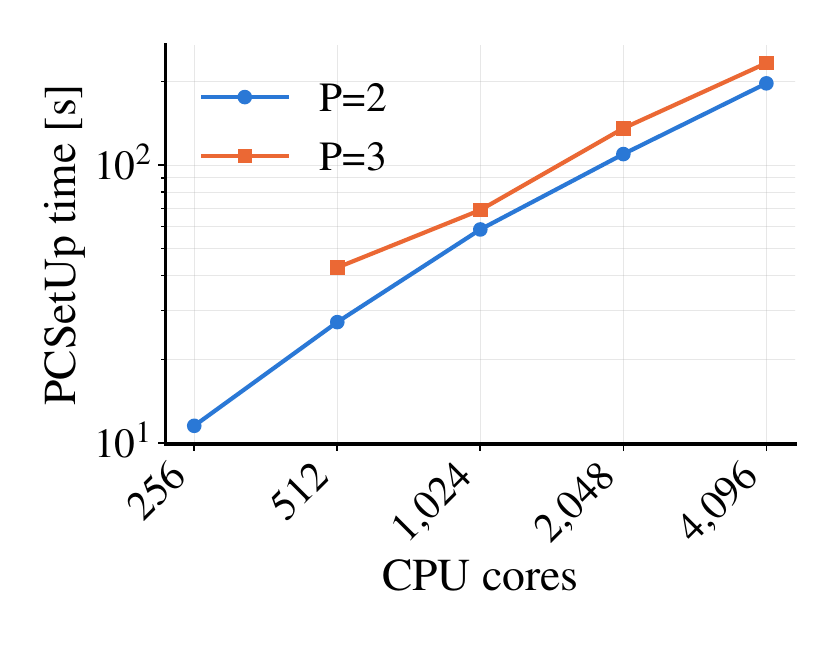}
  \caption{AMG/LOR preconditioner setup time}
  \label{fig:wifw_scaling_pcsetup}
\end{subfigure}

\caption{Strong scaling of the LOR-preconditioned Poisson solve on the IFW-W case.
The outer Krylov solve (\texttt{KSPSolve}) time excludes the one-time
AMG/LOR preconditioner construction (\texttt{PCSetUp}) cost.}
\label{fig:wifw_scaling_time_eff}

\end{figure}

A strong-scaling study of the Poisson pressure system (Eqn.~\eqref{eq:poisson_combined}) is performed on the IFW-W case, with zero Neumann BCs on all boundaries but the outlet, which is specified as zero Dirichlet.
The LOR-preconditioned CG solver is run on ARCHER2 at polynomial orders $P=2$ ($9.6\times10^6$ global DOF) and $P=3$ ($32.2\times10^6$ global DOF), one MPI rank per core.
At $P=3$, 256 cores placed on four or eight nodes exhausted available memory before completing a solve (annotated in Fig.~\ref{fig:wifw_scaling_iterations}); since this occurs regardless of placement, it reflects a genuine memory limit for the $P=3$ problem size at that core count on this mesh, so the $P=3$ sweep begins at 512 cores.
The LOR space is solved with a single AMG V-cycle per outer iteration, using \textit{BoomerAMG} from Hypre via its PETSc interface, with the settings listed in Appendix~\ref{App:petscrc_settings}; this inner-iteration strategy follows the \textit{1iter} configuration validated against the more expensive \textit{tol1e-4} alternative in \cref{sec:LOR_iterative_performance}.

Fig.~\ref{fig:wifw_scaling_iterations} shows that the number of outer CG iterations needed is essentially independent of core count for both polynomial orders, confirming that the conditioning behaviour established in \cref{sec:LOR_conditioning_iterative} carries over to this industrial-scale problem.

Fig.~\ref{fig:wifw_scaling_pcsetup} shows the one-time cost of constructing the AMG (\texttt{PCSetUp}), separate from the Krylov solve (\texttt{KSPSolve}) reported in Fig.~\ref{fig:wifw_scaling_time_eff}.
The setup grows steeply and monotonically with core count, from 12s to 197s for $P=2$ and from 43s to 233s for $P=3$, reflecting the increasing communication cost of building the AMG hierarchy over more MPI ranks.
In an isolated, cold-started benchmark this cost dominates the total wall-clock time, but that is misleading for the intended use case: in a production simulation the pressure operator is constant-coefficient, so the preconditioner is built once and reused over thousands of subsequent timesteps, see \cref{subsec:incns_hx1_note}.
We therefore exclude the setup time from the time-to-solution and efficiency metrics used below, and assess strong scaling using \texttt{KSPSolve} time alone, see Fig.~\ref{fig:wifw_scaling_time_eff}.

Fig.~\ref{fig:wifw_scaling_time_eff} shows that, without setup cost, Krylov solve time decreases essentially monotonically with core count for both polynomial orders.
For $P=2$ it falls from 14.3s at 256 cores to 4.7s at 2048 cores, after which it rises only slightly, to 5.3s at 4096 cores.
For $P=3$ it falls from 40.1s at 512 cores to 16.3s at 4096 cores, with the sweep effectively plateauing beyond 2048 cores.
Parallel efficiency computed on this basis at 2048 cores is 38.2\% for $P=2$ and 59.9\% for $P=3$.
The residual loss of Krylov solve efficiency beyond roughly 1024 cores is consistent with the solve becoming communication-bound as the local subdomain shrinks -- to only $2.3\times10^3$ ($P=2$) and $7.9\times10^3$ ($P=3$) DOF per rank at 4096 cores -- rather than with any loss of algorithmic effectiveness, since the iteration count itself (Fig.~\ref{fig:wifw_scaling_iterations}) is unaffected by core count over the same range.
Together with the amortisation argument above, this shows that once the AMG build cost is paid once and reused over thousands of timesteps, the LOR-preconditioned pressure solve continues to benefit from additional parallelism at least up to the 2048--4096 core range tested here.

%% file: main_computational.tex
\subsection{Incompressible Navier-Stokes solver performance}
\label{subsec:incns_hx1_note}

This section evaluates the performance of the LOR preconditioner in the \IncNS solver of \nekpp, using a production run of the IFW-W case, see Fig.~\ref{fig:testcases_IFW}, at $Re = 2\times10^5$ on ARCHER2.
Three pressure-system configurations, each using an outer conjugate gradient solver, are compared over 1000 timesteps:
\begin{itemize}
    \item \textit{Baseline}: static condensation with Diagonal preconditioner.
    \item \textit{Baseline-full}: full matrix with Diagonal preconditioner.
    \item \textit{LOR}: full matrix with LOR/AMG preconditioner.
\end{itemize}
All three runs use identical mesh, boundary conditions, timestep size, and velocity-solver settings (Appendix~\ref{App:incNS_setup}), and are placed on 2048 MPI ranks, matching the $P=3$ point of the strong-scaling sweep in \cref{subsec:wifw_strong_scaling}.\footnote{Each configuration is reported from a single such run rather than an ensemble due to high computational costs.}

\begin{figure}[H]
\centering
\begin{subfigure}{.48\textwidth}
  \centering
  \includegraphics[width=\linewidth]{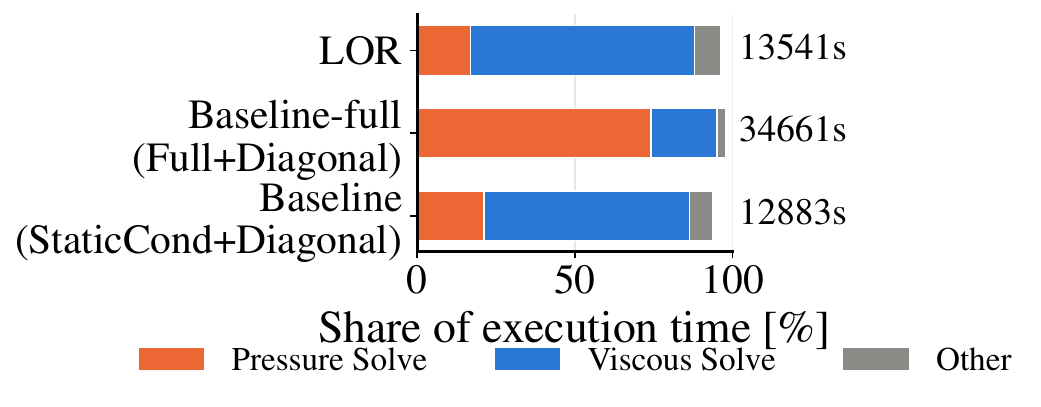}
  \caption{Execution time breakdown}
  \label{subfig:IFW_W_timesplit}
\end{subfigure}%
\begin{subfigure}{.48\textwidth}
  \centering
  \includegraphics[width=\linewidth]{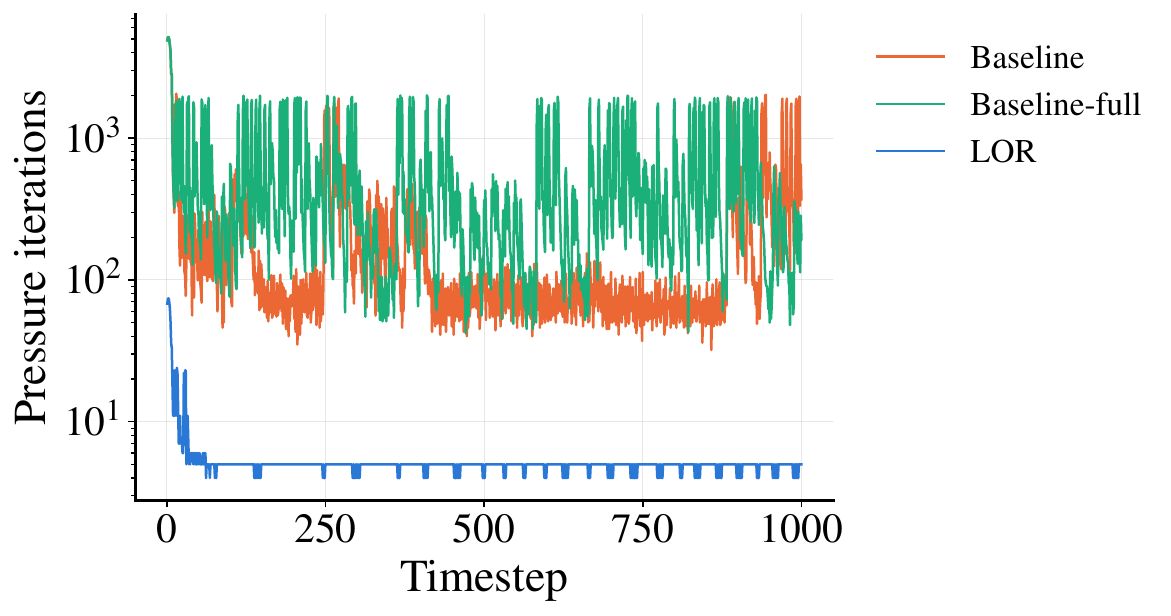}
  \caption{Pressure iterations per step}
  \label{subfig:IFW_W_pressure_iterations}
\end{subfigure}
\caption{Production \IncNS run on the IFW-W case matching Table~\ref{tab:wifw_incns_production}.
(a) Share of total execution time spent in the pressure solve, viscous solve, and all other operators (advection, projections, boundary-condition and forcing assembly), with total wall-clock time annotated.
(b) Outer pressure-solve iteration count at every one of the 1000 timesteps, showing both the mean reduction and the much lower step-to-step variance of the LOR preconditioner relative to Diagonal.}
\label{fig:IFW_W_production}
\end{figure}

Table~\ref{tab:wifw_incns_production} and Fig.~\ref{fig:IFW_W_production} show that LOR reduces the pressure-solve wall-clock time by $11.2\times$ relative to Baseline-full. The pressure solve accounts for 74\% of the total execution time with Baseline-full, compared with 17\% with LOR. This reduction is driven by the decrease in the mean outer iteration count from 559.6 to 5.5, corresponding to a $101\times$ reduction, together with a substantially smaller step-to-step variation, as shown in Fig.~\ref{subfig:IFW_W_pressure_iterations}. Baseline and Baseline-full vary between approximately 30 and 5200 iterations as the local flow evolves, whereas the LOR iteration count remains within a factor of approximately $19$ of its minimum. Over the extended 1000-step run, the mean LOR iteration count decreases to 4.9 over the final 500 steps and 4.8 over the final 100 steps. These values provide a more representative estimate of the steady production regime after the initial transient.

\begin{table}[H]
\centering
\small
\caption{Production \IncNS run on the IFW-W case using 2048 MPI ranks, polynomial order $P=3$ and 1000 timesteps.
Percentages are the share of total execution time.
Pressure iterations are the mean outer CG iteration count for the $p$-solve with the observed min--max range in parenthesis, velocity iterations are the mean outer CG iteration count for $u,v,w$ combined.}
\label{tab:wifw_incns_production}
\begin{tabular}{lccc}
\hline
 & \textbf{Baseline} & \textbf{Baseline-full} & \textbf{LOR} \\
 & (StaticCond+Diag) & (Full+Diag) & (Full+LOR) \\ \hline
Total Time [s]              & 12{,}883      & 34{,}661      & 13{,}541 \\
Pressure Solve [s] (\%)     & 2727 (21.1\%) & 25{,}699 (74.1\%) & 2289 (16.9\%) \\
Viscous Solve [s] (\%)      & 8414 (65.2\%) & 7273 (21.0\%) & 9640 (71.1\%) \\
Pressure Iters (mean, range) & 235.3 (32--5156) & 559.6 (42--5164) & 5.5 (4--74) \\
Velocity Iters (mean, $u,v,w$) & 69.7 & 69.7 & 69.7 \\ \hline
\end{tabular}
\end{table}

Relative to Baseline, LOR reduces the mean outer iteration count from 235.3 to 5.5, a factor of \(42\), but reduces the pressure-solve wall-clock time from \(2727\,\mathrm{s}\) to \(2289\,\mathrm{s}\), a factor of only \(1.2\). Each LOR iteration applies a full AMG V-cycle to the large, non-statically-condensed LOR system and costs \(0.41\,\mathrm{s}\). Each Baseline iteration uses a diagonal-preconditioned matrix-vector product on the smaller, statically-condensed system and costs \(0.012\,\mathrm{s}\). This \(34\times\) difference in per-iteration cost offsets most of LOR's iteration-count advantage.

Static condensation and LOR/AMG address different computational costs in this regime. LOR provides its largest speedup over the full system for which it was designed. Against the inexpensive, statically-condensed Baseline, it yields a modest reduction in wall-clock time while requiring far fewer pressure iterations, as shown in Fig.~\ref{subfig:IFW_W_pressure_iterations}.


All three configurations require an average of 69.7 velocity iterations, indicating that the pressure-preconditioner comparison does not alter velocity-solve convergence. The variation in viscous-solve time in Table~\ref{tab:wifw_incns_production} is not investigated further here; each configuration is represented by a single production run, as noted above.

These production results corroborate the amortisation argument of \cref{subsec:wifw_strong_scaling}: the one-time AMG build cost is essentially unchanged from the isolated, cold-started $P=3$ benchmark at the same core count and DOF (130s here vs.\ 136s), confirming that this cost depends on the operator and its parallel decomposition alone.
What differs is the outer iteration count: thanks to initial-guess extrapolation across timesteps following the work in \cite{FISCHER1998193}, only the first few of the 1000 solves approach the cold-start cost of 54 iterations, against a production mean of 5.5.
The resulting effective per-step LOR cost ($2289\text{s}/1000 = 2.3$s) is far below the 381s a naive extrapolation from the single cold-started benchmark would suggest.

%% file: appendix.tex
\section{Linear solvers in \nekpp}
\label{App:linear_solver_nekpp}
Large-scale \nekpp simulations utilize iterative Krylov subspace solvers targeting the Helmholtz system until the $k$-th iteration residual, $\hat{\mathbf{r}}_{k} = \mathbf{f} - \mathbf{H}\hat{\mathbf{u}}_{k}$, satisfies a prescribed tolerance $\tau$. While relative criteria ($\|\hat{\mathbf{r}}_{k}\|_{2} < \tau \|\hat{\mathbf{r}}_{0}\|_{2}$) are common, their sensitivity to problem configuration often necessitates geometry-specific tuning. To ensure consistent convergence behavior across diverse scales and geometries, \nekpp adopts an absolute $L_2$ termination criterion $\|\hat{\mathbf{r}}_{k}\|_{2} < \tau$. Native implementations include conjugate gradient (CG) and Generalized Minimal Residual (GMRES) methods \cite{Saad2003}, while a broader range of linear solvers and preconditioners is available through the PETSc interface \cite{petsc-web-page}. This interface provides access to algebraic multigrid (AMG) \cite{stuben2001review} libraries, such as BoomerAMG \cite{yang2002boomeramg}, which are particularly advantageous for unstructured meshes and complex cases where geometric information is unavailable.

\section{Software}
\label{App:software}
All numerical experiments using the LOR preconditioner were performed with the master branch of \nekpp. The corresponding implementation is provided in Nektar++ merge request \href{https://gitlab.nektar.info/nektar/nektar/-/merge_requests/1888}{!1888}.

In the HX1 cluster, the code is linked against an independently compiled version of PETSc v3.22.0, configured with HYPRE v2.31.0 to enable the use of BoomerAMG routines for algebraic multigrid solves. The software stack is built with the \texttt{intel/2022a} toolchain, which provides access to optimised BLAS, LAPACK, and MPI implementations. In ARCHER2, PETSc v3.19.3 is used, bundled as a Third-Party library within the \nekpp build system. This version links to HYPRE v2.28.0. The code is compiled using the GNU Compiler Collection (GCC) under the \texttt{PrgEnv-gnu} environment.

\begin{table}[H]
\centering
\caption{Summary of software environment used on HX1 and ARCHER2 clusters.}
\label{tab:software_env}
\begin{tabular}{|l|l|l|}
\hline
\textbf{Component}        & \textbf{HX1} & \textbf{ARCHER2} \\ \hline
\nekpp                 & Development branch with LOR MR & Same \\ \hline
PETSc                    & v3.22.0 (independently compiled) & v3.19.3 (ThirdParty) \\ \hline
HYPRE                    & v2.31.0                          & v2.28.0                     \\ \hline
Toolchain                & \texttt{intel/2022a}             & \texttt{PrgEnv-gnu}         \\ \hline
MPI Library              & Intel MPI                        & Cray MPICH                  \\ \hline
\end{tabular}
\end{table}

\section{Hardware}
\label{App:hardware}
The studies presented in this work were carried out on the following high-performance computing (HPC) platforms.

\textbf{HX1}: \href{https://icl-rcs-user-guide.readthedocs.io/en/latest/hpc/cluster-specification/#hx1}{HX1} (or Hex) is the production``capability" HPC cluster hosted by the Imperial College Research Computing Service. Each compute node comprises a Lenovo SD630v2 server equipped with two Intel Xeon Platinum 8358 (Ice Lake) processors, operating at 2.60 GHz with 32 cores per socket, yielding a total of 64 cores and 512 GB of RAM per node. Nodes are interconnected via NVIDIA ConnectX-6 HDR200 InfiniBand, offering a peak interconnect bandwidth of 200 Gbit/s. 

\textbf{ARCHER2}: \href{https://www.archer2.ac.uk/about/hardware.html}{ARCHER2}, the UK national supercomputing service, consists of compute nodes featuring dual AMD EPYC 7742 processors. Each processor provides 64 cores at 2.25 GHz, and all cores are fully utilised during the runs. ARCHER2 nodes are configured with 256 GB of RAM and communicate via the HPE Cray Slingshot interconnect, delivering a bidirectional bandwidth of 2 × 100 Gbps per node. 

\section{PETSc settings}
\label{App:petscrc_settings}
\begin{itemize}[noitemsep]
    \item 1 V-cycle for AMG, or until a tolerance, if specified
    \item HMIS coarsening
    \item Strong threshold = 0.70 for 3D cases
    \item Coarse grid size = 10 to control the levels of AMG cycle
    \item Extended+i interpolation with a Pmax = 4 and truncation of 0.3
    \item L1Scaled SOR/Jacobi smoother with two sweeps going up and down
\end{itemize}

\section{\IncNS solver simulation setup}
\label{App:incNS_setup}
The simulation is setup in accordance with the methodology detailed in \cite{khurana2025industrialization} The computational domain emulates the geometry inside a wind tunnel. The freestream air enters the domain and is enforced as a uniform Dirichlet boundary condition at the inlet. The ceiling and the side walls are treated as slip walls. A stable high-order form of the zero Neumann boundary condition at the outlet is used for the velocity, and the pressure is fixed to the ambient pressure.

\begin{itemize}[noitemsep]
    \item Semi-implicit time stepping is used \cite{karniadakis2005spectral}
    \item Taylor-Hood approximation: The polynomial order for the velocity is \textit{P+1}, pressure is \textit{P}. The simulations are run at a polynomial order of \textit{P=4,3}.
    \item A second-order time-integration scheme with the maximum constant stable time-step of $5\times10^{-6}\,\mathrm{s}$. 
    \item Initialisation from a velocity field obtained from a restart solution, and the pressure is 0.
\end{itemize}
The pressure system contains $1.12\times10^8$ local and $3.22\times10^7$ global degrees of freedom, while each velocity component contains $2.04\times10^8$ local and $7.62\times10^7$ global degrees of freedom. Including all three velocity components, the complete system contains $7.24\times10^8$ local and $2.608\times10^8$ global degrees of freedom.

\raggedbottom
\section{Construction of the low-order refined mesh}
\label{app:lor_mesh_construction}

This appendix details the construction of the LOR mesh introduced in Section~\ref{subsec:lor_space_construct}. For a degree-$P$ parent expansion, the $P+1$ points on each tensor-product edge define $n_{\mathrm{split}}=P$ refined segments. The construction preserves the existing global identifiers and orientations on shared edges and faces, creates the required interior entities, and updates the coefficient maps used to transfer data between the HO and LOR representations. Algorithm~\ref{alg:linmesh-setup-2d} describes the common two-dimensional face construction. Algorithms~\ref{alg:linmesh-setup-tet} and~\ref{alg:linmesh-setup-prism} then use these face connectivities to construct tetrahedral and prismatic subelements, respectively.

\input{algos/build_LinMeshSetUp2DGeom}

\input{algos/build_LinMeshSetUpTetGeom}

\input{algos/build_LinMeshSetUpPrismGeom}

%% file: algos/build_LinMeshSetUp2DGeom.tex
\begin{algorithm}[H]
\caption{Two-dimensional LOR mesh construction}
\label{alg:linmesh-setup-2d}
\begin{algorithmic}[1]
\REQUIRE HO mesh graph $\mathcal{M}_H$, LOR mesh graph $\mathcal{M}_L$, subdivision count $n_{\mathrm{split}}$, face--edge offset map $\mathcal{O}_{FE}$, two-dimensional coefficient map $\mathcal{C}_{2D}$, vertex offset $v_{\mathrm{off}}$, prescribed nodal distribution, and refined-cell type
\ENSURE Updated vertex, edge, and face connectivity in $\mathcal{M}_L$, together with the updated coefficient map $\mathcal{C}_{2D}$

\STATE \textbf{Stage 1: initialise data structures}
\STATE Allocate coordinate arrays $X$, $Y$, and $Z$, and local connectivity arrays for vertices, horizontal edges, vertical edges, and diagonal edges
\STATE Determine the maximum existing global vertex and edge identifiers used to offset newly created entities
\STATE Compute the edge-index offset $n_{\mathrm{split}}(3n_{\mathrm{split}}-2)$

\STATE \textbf{Stage 2: process triangular faces}
\FOR{each triangular face $\mathcal{F}_t$ in $\mathcal{M}_H$}
    \STATE Evaluate the physical coordinates at the points of the prescribed nodal distribution
    \STATE Retrieve the global edge identifiers and their orientations
    \STATE Assign boundary vertex identifiers using $n_{\mathrm{split}}$ subdivisions on each edge
    \STATE Create and index the interior vertices
    \STATE Construct the subdivided boundary, interior, and diagonal edges
    \STATE Insert the resulting linear triangular subelements into $\mathcal{M}_L$
\ENDFOR

\STATE \textbf{Stage 3: process quadrilateral faces}
\FOR{each quadrilateral face $\mathcal{F}_q$ in $\mathcal{M}_H$}
    \STATE Evaluate the physical coordinates on the selected tensor-product nodal grid
    \STATE Retrieve the global edge identifiers, orientations, and local ordering
    \STATE Assign boundary and interior vertex identifiers
    \STATE Construct the horizontal and vertical edges of the refined grid
    \IF{the prescribed refined-cell type is triangular}
        \STATE Add one diagonal edge to each refined quadrilateral cell
        \STATE Insert two triangular subelements per refined cell
    \ELSE
        \STATE Insert one quadrilateral subelement per refined cell
    \ENDIF
\ENDFOR

\STATE Update $\mathcal{O}_{FE}$ and $\mathcal{C}_{2D}$ with the new local-to-global connectivity
\RETURN $\mathcal{M}_L$ and $\mathcal{C}_{2D}$
\end{algorithmic}
\end{algorithm}

%% file: algos/build_LinMeshSetUpTetGeom.tex
\begin{algorithm}[H]
\caption{Tetrahedral LOR mesh construction}
\label{alg:linmesh-setup-tet}
\begin{algorithmic}[1]
\REQUIRE HO mesh graph $\mathcal{M}_H$, LOR mesh graph $\mathcal{M}_L$, subdivision count $n_{\mathrm{split}}$, face--edge offset map $\mathcal{O}_{FE}$, and three-dimensional coefficient map $\mathcal{C}_{3D}$
\ENSURE Updated vertex, edge, face, and volume connectivity for the linear tetrahedral subelements in $\mathcal{M}_L$

\STATE \textbf{Stage 1: initialise data structures}
\STATE Allocate local indexing arrays for vertices, edges, and faces
\STATE Determine the maximum existing global identifiers used to offset newly created entities

\FOR{each tetrahedron $\mathcal{T}_e$ in $\mathcal{M}_H$}
    \STATE \textbf{Stage 2: determine parent-element connectivity}
    \STATE Retrieve the vertex coordinates and local face--edge orientations of $\mathcal{T}_e$
    \STATE Initialise the local indexing of subdivided vertices and edges
    \FOR{each face $F_m$ of $\mathcal{T}_e$}
        \STATE Assign its face--edge connectivity using $\mathcal{O}_{FE}$
        \STATE Reuse the vertex and edge connectivity already stored for $F_m$ in $\mathcal{M}_L$
    \ENDFOR

    \STATE \textbf{Stage 3: compute subdivision-point coordinates}
    \STATE Evaluate the physical-coordinate mapping at the selected simplex nodal points on the edges, faces, and interior of $\mathcal{T}_e$

    \STATE \textbf{Stage 4: construct the tetrahedral subelements}
    \STATE Create and index the new interior vertices
    \STATE Construct the linear edges using the local connectivity and $\mathcal{O}_{FE}$
    \STATE Construct consistently oriented triangular faces from the linear edges
    \STATE Assemble the tetrahedral subelements from the triangular faces

    \STATE \textbf{Stage 5: update the global LOR mesh}
    \STATE Insert the new vertices, edges, faces, and tetrahedra into $\mathcal{M}_L$
    \STATE Update $\mathcal{C}_{3D}$ and the geometric connectivity
\ENDFOR

\RETURN $\mathcal{M}_L$ and $\mathcal{C}_{3D}$
\end{algorithmic}
\end{algorithm}

%% file: algos/build_LinMeshSetUpPrismGeom.tex
\begin{algorithm}[H]
\caption{Prismatic LOR mesh construction}
\label{alg:linmesh-setup-prism}
\begin{algorithmic}[1]
\REQUIRE HO mesh graph $\mathcal{M}_H$, LOR mesh graph $\mathcal{M}_L$, subdivision count $n_{\mathrm{split}}$, face--edge offset map $\mathcal{O}_{FE}$, and three-dimensional coefficient map $\mathcal{C}_{3D}$
\ENSURE Updated vertex, edge, face, and volume connectivity for the linear prismatic subelements in $\mathcal{M}_L$

\FOR{each prism $\mathcal{P}_e$ in $\mathcal{M}_H$}
    \STATE \textbf{Stage 1: determine parent-element connectivity}
    \STATE Retrieve the vertex, edge, and face definitions of $\mathcal{P}_e$
    \STATE Determine the orientation and alignment of its two triangular faces
    \IF{the triangular faces are not consistently aligned}
        \STATE Report that the prism cannot be subdivided with the required orientation
    \ENDIF

    \STATE \textbf{Stage 2: subdivide edges and faces}
    \FOR{each edge $E_i$ of $\mathcal{P}_e$}
        \STATE Subdivide $E_i$ into $n_{\mathrm{split}}$ linear segments using the prescribed nodal distribution
    \ENDFOR
    \FOR{each face $F_j$ of $\mathcal{P}_e$}
        \IF{$F_j$ is triangular}
            \STATE Construct $n_{\mathrm{split}}^2$ triangular subfaces
        \ELSE
            \STATE Construct $n_{\mathrm{split}}^2$ quadrilateral subfaces
        \ENDIF
    \ENDFOR

    \STATE \textbf{Stage 3: construct the prismatic subelements}
    \STATE Connect corresponding triangular subdivisions through the one-dimensional refined direction
    \STATE Generate the linear prisms with orientations inherited from $\mathcal{P}_e$
    \STATE Record the new vertex and edge offsets in $\mathcal{O}_{FE}$

    \STATE \textbf{Stage 4: update the global LOR mesh}
    \STATE Insert the new vertices, edges, faces, and prisms into $\mathcal{M}_L$
    \STATE Update $\mathcal{C}_{3D}$ and the geometric connectivity
\ENDFOR

\RETURN $\mathcal{M}_L$ and $\mathcal{C}_{3D}$
\end{algorithmic}
\end{algorithm}